\documentclass{article}

\usepackage{amsmath, amscd, amssymb, graphics, xypic, mathrsfs}

\newcommand{\setof}[1]{\{ #1 \}}

\newcommand{\tensor}{\otimes}

\newcommand{\colim}{\operatorname{colim}}
\newcommand{\Spec}{\operatorname{Spec}}
\newcommand{\isomto}{{\stackrel{\sim}{\longrightarrow}}}

\renewcommand{\O}{{\mathcal O}}

\renewcommand{\L}{{\mathcal L}}

\newcommand{\F}{{\mathcal F}}

\newcommand{\ho}[1]{{\mathcal H}({#1})}
\newcommand{\hop}[1]{{\mathcal H}_{\cdot}({#1})}
\newcommand{\dmeff}[1]{{\bf DM}^{eff,-}_{Nis}(k,{#1})}

\newcommand{\cplx}{{\mathbb C}}

\newcommand{\Q}{{\mathbb Q}}
\newcommand{\Z}{{\mathbb Z}}
\newcommand{\aone}{{\mathbb A}^1}
\newcommand{\pone}{{\mathbb P}^1}

\newcommand{\gm}{{{\mathbb G}_{\bf m}}}

\newcommand{\et}{\text{\'et}}

\newcommand{\Sm}{{\mathcal Sm}}

\newcommand{\Shv}{{\mathcal Shv}}

\newcommand{\Spc}{{\mathcal Spc}}
\newcommand{\Fq}{{{\mathbb F}_q}}
\newcommand{\Ox}{{\mathcal O}_{X,x}}
\newcommand{\Ohat}{\hat{{\mathcal O}_{x}}}
\newcommand{\Khat}{\hat{{\mathcal K}_{x}}}

\newcommand{\umu}{{{\mu}}}
\newcommand{\cochar}[1]{{\mathbf{X}_*(#1)}}
\newcommand{\chara}[1]{{\mathbf{X}^*(#1)}}
\newcommand{\underlinemu}{{\mu}}
\renewcommand{\F}{{\mathcal F}}

\newcommand{\PP}{{\mathbb P}}

\newcommand{\calc}{{\mbox{$\mathcal C$}}}
\newcommand{\calm}{{\mbox{$\mathcal M$}}}
\newcommand{\mnd}{{{\mathcal M}(n,d)}}

\newcommand{\mndc}{{{\mathcal M}^C(n,d)}}

\newcommand{\divnd}{{{\rm Div}^{n,d}_{C/k}}}
\newcommand{\divt}{{{\rm Div}^{n,\tilde{d}}_{C/k}}}

\newtheorem{thm}{Theorem}[section]
\newtheorem{lem}[thm]{Lemma}
\newtheorem{cor}[thm]{Corollary}
\newtheorem{prop}[thm]{Proposition}

\newtheorem{defn}[thm]{Definition}

\newtheorem{rem}[thm]{Remark}
\newtheorem{ex}[thm]{Example}

\begin{document}

\title{Yang-Mills theory and Tamagawa Numbers: \\ \scriptsize{the
fascination of unexpected links in mathematics}}

\author{Aravind Asok, Brent Doran and Frances Kirwan}

%\classno{14H60  (14F42 14L24)}

%\begin{document}
\maketitle

\section{Introduction and overview}
This article is an expanded version of the third author's
presidential address to the London Mathematical Society in November
2005, which was based on work with the other two authors over the
preceding year.  The title was inspired by Michael Atiyah's
presidential address delivered in 1976 \cite{mfa}
 in which he said
\begin{quote}
``The aspect of mathematics which fascinates me most is the rich
interaction between its different branches, the unexpected links,
the surprises."
\end{quote}

The unexpected link which is the topic of this article was remarked
on by Atiyah himself and his collaborator Raoul Bott in their
fundamental 1983 paper \cite{AB} on the Yang-Mills equations over
Riemann surfaces.  In this paper Atiyah and Bott used ideas coming
from Yang-Mills theory and equivariant Morse theory to derive
inductive formulae for the Betti numbers of the moduli spaces $\mnd$
of stable vector bundles of rank $n$ and degree $d$ over a fixed
compact Riemann surface $C$ of genus $g \geq 2$, when $n$ and $d$
are coprime. (We will assume throughout this
introduction that $n$ and $d$ are coprime integers with $n>0$.)
Equivalent formulae had been obtained earlier by Harder and
Narasimhan \cite{HN} and Desale and Ramanan \cite{DR} using
arithmetic techniques and the Weil conjectures.  In the latter
approach a crucial ingredient was the fact, proved by Weil, that
the volume of a certain locally symmetric space attached to $SL_n$
with respect to a canonical measure
-- an invariant known as the Tamagawa number of $SL_n$ -- is $1$.

Atiyah and Bott observed that although the two methods came from
very different branches of mathematics, namely arithmetic and
physics, there was a formal correspondence between them, with the
Tamagawa number of $SL_n$ (or equivalently the function field
analogue of the Siegel formula) playing, roughly speaking, the r\^{o}le
of the cohomology of the classifying space of the gauge group in the
Atiyah-Bott approach.  They asked for a deeper understanding of this
observation and in particular for a geometric explanation, exploiting
the analogy with equivariant cohomology, of the fact that the
Tamagawa number of $SL_n$ is $1$.  Contributions since then towards
such understanding have included work by Bifet, Ghione and
Letizia \cite{Bifet,BGL}, providing another inductive procedure for
calculating the Betti numbers of the moduli spaces which is in some
sense intermediate between the arithmetic approach and the
Yang-Mills approach, and more recently, work by Teleman, Behrend,
Dhillon and others on the moduli stack of bundles over $C$
\cite{BD05b,BD05a,Dh06,Tel98}.

The formal correspondence observed by Atiyah and Bott between the
inductive formulae obtained from the Yang-Mills and arithmetic
points of view arises because both depend on stratifications of
spaces whose points represent vector bundles over $C$, with the
stratification induced by the ``Harder-Narasimhan types" of the
bundles.  In each case it is the open stratum, corresponding to
semistable bundles, which needs to be understood, and the inductive
calculation comes from combining a simple description of the other
strata (in terms of semistable strata for the corresponding problem
with strictly smaller values of $n$ and varying values of $d$)
together with knowledge of the whole space.  In the arithmetic
approach the space to be stratified is a coset space of the group
$SL_n(\mathbb{A}_K)$ associated with the ad\`ele ring
${\mathbb{A}}_K$ of the function field $K$ of a nonsingular
projective curve over a finite field.   The Tamagawa measure of the
whole space is the (countably infinite) sum of the measures of the
strata. In the Yang-Mills approach Atiyah and Bott stratify the
infinite-dimensional affine space $\mathcal{A}$ of unitary
connections on a fixed $C^{\infty}$ bundle of rank $n$ and degree
$d$ over $C$, and they show that the stratification is \lq\lq
equivariantly perfect" with respect to the gauge group, so that the
equivariant cohomology algebra of $\mathcal{A}$ is isomorphic as a
vector space to the sum over the strata of their equivariant
cohomology algebras, shifted in degree by their real codimensions.

The cohomology of the moduli spaces $\mnd$ of bundles over a compact
Riemann surface has, of course, more structure than is revealed by
its Betti numbers. It has been an object of study for several
decades, first when $n=2$ and then for general $n$, from many points
of view including mathematical physics \cite{AB,Witten}, matrix
divisors \cite{BGL,Bifet}, extended moduli spaces
\cite{JW94,HJ94,J}, group valued moment maps and quasi-Hamiltonian
reductions \cite{AMM98,AMW00,AMW01,AMW02,MW99,Mei05} and others
\cite{Dh06,KirMaps,Kir5,JKInt,JKKW2}, though the ideas of Atiyah and
Bott have been fundamental for much of this progress.

Following a long line of work initiated by Grothendieck in the
1960s, Bloch and Voevodsky have developed sophisticated versions of
cohomology (motivic cohomology, see \cite{Bloch,VTriCat,MVW}) for
algebraic varieties defined over arbitrary fields.  Voevodsky and
Morel had a broader vision (see \cite{MV}): they sought to build a
homotopy theory for schemes over a field $k$ where the affine line
plays a role analogous to that played by the unit interval in
ordinary homotopy theory.  In the resulting homotopy category,
called the $\aone$ or motivic homotopy category, there are two
different analogues of the circle: the simplicial circle which is
obtained from the affine line by glueing together $0$ and $1$
(represented by the affine nodal cubic curve) and the Tate circle
$\aone - \{0\}$.  Corresponding to these two analogues of the
circle, a pair of integers index the motivic cohomology groups of a
smooth scheme $X$ over a field $k$, reflecting the fact that these groups
form a bigraded ring. Motivic cohomology groups for smooth schemes
over $k$ have many properties analogous to ordinary singular
cohomology, including an appropriate form of homotopy invariance,
Mayer-Vietoris and Gysin long exact sequences. However, motivic
cohomology encodes far more information; for example, in contrast
with ordinary singular cohomology, the motivic cohomology of a point
$\Spec k$ is quite large, with its degree $(p,p)$ part being
isomorphic to the $p$-th Milnor $K$-group of the field $k$ (see
\cite{MVW} Chapter 4). Furthermore, if $X$ is a smooth $k$-scheme,
the motivic cohomology groups $H^{2p,p}(X,\Z)$ are isomorphic to the
Chow groups of codimension $p$ cycles on $X$ modulo rational
equivalence.

The aim of this article is threefold: to announce results, produced
in the setting of $\mathbb{A}^1$-homotopy theory, on the motivic
cohomology (and hence its many realizations) of quotients in the
sense of Mumford's geometric invariant theory or GIT \cite{ADK1},
generalizing what was known for singular cohomology of GIT
quotients;
 to explain how methods used over the last three decades to study the
singular cohomology of the moduli spaces $\mnd$  can thence be
adapted to study the motivic cohomology of moduli spaces of bundles
on a smooth projective curve $C$ over an algebraically closed field
$k$;
 and at the same time to extend our understanding of the link
between Yang-Mills theory and Tamagawa numbers remarked on by
Atiyah and Bott, by re-examining some of their essentially homotopy
theoretic considerations in a more algebraic modern light.

This study is based on an adaptation to the setting of motivic
cohomology of the inductive methods obtained in the third author's
thesis \cite{Kir1} for calculating the Betti numbers of a GIT quotient
\cite{GIT} of a nonsingular complex projective variety $X$ by a
linear action of a complex reductive group $G$.  These methods were
themselves inspired by the work of Atiyah and Bott \cite{AB} and involve applying equivariant
Morse theory to the norm-square of an appropriate moment map, which
is the analogue of the Yang-Mills functional in \cite{AB}.  The
associated stratification of $X$ has an alternative purely algebraic
description, independent of Morse theory, which is valid much more
generally than just over the complex numbers $\cplx$: the semistable
points of $X$ (in the sense of GIT) form an open stratum, and the
other strata can be described inductively in terms of the semistable
points of nonsingular projective subvarieties of $X$ under
appropriately linearized actions of reductive subgroups of $G$. In
the finite-dimensional algebro-geometric setting, most of the
results of \cite{Kir1} can be adapted to apply to motivic cohomology
\cite{ADK1}. These results can then be used to study the motivic
cohomology of moduli spaces of vector bundles over a curve $C$.

Much beautiful work has been done recently on the moduli of bundles
over a curve which can be used to extract similar information
\cite{dBmotive,dBrk,BD05b,Dh06}. Our aim here is to derive
this information as an application of a general theory of cohomology
of GIT quotients, and to offer some insights gleaned from returning
to a mapping space perspective on the moduli space of bundles.

\subsubsection*{Outline}

In \S 2 of this paper we describe the two equivalent inductive
procedures for calculating the Betti numbers of $\mnd$ provided by
the arithmetic and Yang-Mills approaches. The Yang-Mills approach of
Atiyah and Bott \cite{AB}, which uses equivariant Morse theory, has
as its basic ingredient a simple description given in \cite{AB} of
the cohomology of the classifying space of the gauge group, while
the arithmetic approach of \cite{HN,DR} uses Tamagawa measures and
reduces to the function field version of the Siegel formula, or
equivalently to the fact that the Tamagawa number of $SL_n$ is $1$.

These two approaches rely on what can be regarded as
infinite-dimensional quotient constructions of the moduli spaces
$\mnd$.  A third approach, closely related to that used by Bifet,
Ghione and Letizia in \cite{BGL}, involves the construction of the
moduli spaces $\mnd$ as finite-dimensional GIT quotients which can
be regarded as finite-dimensional approximations to the Yang-Mills
construction. This involves studying spaces of maps into
Grassmannians and matrix divisors, and provides an inductive
calculation of the Betti numbers of the moduli spaces in terms of
the cohomology of symmetric powers of the curve $C$ (this cohomology
is well known \cite{MacDonald}).  Here the inductive formulae can be
obtained from equivariant Morse theory as in the Yang-Mills
approach, but in this finite-dimensional setting the Weil
conjectures provide an alternative derivation by using the algebraic
description of the Morse strata and counting points on associated
varieties defined over finite fields.

In \S 3 we move into the motivic world; our aim is to show that the
finite dimensional GIT construction of $\mnd$ used in the third
approach to calculate Betti numbers can also be used to study the
motivic cohomology of $\mnd$.  The first step is to make sense of
equivariant motivic cohomology; this is done using a straightforward
modification of the Borel construction in topology using work of
Morel and Voevodsky. In \S 4 we describe how to adapt the methods of
\cite{Kir1} on the cohomology of GIT quotients to the setting of
motivic cohomology, and in \S 5  we apply these to obtain an
inductive description of the motivic cohomology groups of $\mnd$ in
terms of the motivic cohomology of products of symmetric powers of
the curve $C$.

In the final section \S 6 we leave motivic cohomology and return to
the original question of the relationship between Yang-Mills theory
and Tamagawa measures. We observe that the inductive procedure for
calculating the Betti numbers of $\mnd$ which uses
finite-dimensional approximations to the Yang-Mills picture, via
maps into Grassmannians and matrix divisors, is directly equivalent
to the Yang-Mills approach through a generalization of Segal's
theorem on the topology of spaces of rational functions
\cite{Segal}. This theorem makes precise the sense in which the
Yang-Mills approach is an infinite-dimensional limit of the
corresponding procedure using maps into Grassmannians. In a similar
way the arithmetic methods used by Harder and Narasimhan to provide
an inductive calculation of the Betti numbers of $\mnd$ can be
regarded as an infinite-dimensional limit of the alternative
finite-dimensional approach which involves counting points on
associated varieties over finite fields. The Weil conjectures then
provide the final link in the chain connecting the Yang-Mills and
arithmetic viewpoints.

\subsubsection*{Acknowledgements}
The authors would like to thank Tamas Hausel, Michael McQuillan, James Parson,
Burt Totaro and
Akshay Venkatesh for helpful discussions and comments.

\section{A review of the classical constructions}
The moduli spaces of stable vector bundles on a smooth projective
curve $C$ over a field $k$ have different incarnations depending on the
field $k$.  Over a field $k$ with sufficiently many elements,
these moduli spaces can be constructed by means of geometric
invariant theory.  If $k = \cplx$ and we view the complex points of
$C$ as a compact Riemann surface, we have a differential
geometric construction using an interpretation of stable vector
bundles in terms of connections. If $k = \Fq$, there is an
interpretation of such vector bundles in terms of {\em ad\`eles}. We
review the last two of these constructions in this section.  We will
discuss the GIT construction in \S5.

Historically, the computation of the cohomology of these moduli
spaces was first achieved using number theory and then using
Yang-Mills theory, while algebraic geometry provides a logical
bridge between these two contexts. In this section we will review
the differential geometric and the arithmetic approaches to
computing the cohomology of the moduli spaces.  Throughout this
section, $C$ will denote a compact Riemann surface or a smooth
projective algebraic curve defined over a field $k$. We abuse
notation in this way to emphasize the interchangeability of the
different points of view, the object under consideration being clear
from its context.

\subsection{Yang-Mills and Riemann surfaces}
Let $C$ be a compact Riemann surface of genus $g$. In order to avoid
having to consider special cases, we will assume throughout that $g$
is at least 2.  By convention, holomorphic vector bundles on $C$
will be denoted by calligraphic letters ${\mathcal E},\F, \ldots$
and the underlying $C^{\infty}$-complex vector bundles will be
denoted by Roman letters $E,F,\ldots$.  Furthermore, all bundles
will be assumed to be $C^{\infty}$-bundles; thus the term complex
vector bundle should be read as $C^{\infty}$-complex vector bundle.
Given a complex vector bundle $E$, the bundle of frames of $E$ is a
$GL_n(\cplx)$-principal bundle.  This construction defines a
bijection between the set of  complex vector bundles of rank $n$ on
$C$ and the set of $GL_n(\cplx)$-principal bundles (with inverse
given by forming the vector bundle associated with the standard
$n$-dimensional representation of $GL_n$).  By abuse of terminology,
we will use the same notation $E$ for the $GL_n(\cplx)$-principal
bundle associated with a complex vector bundle $E$.

Topologically, complex vector bundles of fixed rank $n$ on $C$ are
classified by homotopy classes of maps to the space $BGL_n$, which
is (homotopy equivalent to) the Grassmannian of $n$-dimensional
quotients of an infinite dimensional complex vector space.  The
degree or first Chern class $d \in H^2(C,\Z)$, which can be
identified with $\Z$ by pairing with the fundamental class, is also
a topological invariant and the rational number $\mu(E)=d/n$ is
called the {\em slope} of $E$.  We will also write $\mu({\mathcal
E})$ for the slope of the complex vector bundle $E$ underlying a
holomorphic vector bundle ${\mathcal E}$.

If ${\mathcal E}$ is a holomorphic vector bundle, we call ${\mathcal
E}$ {\em stable} (respectively {\em semistable}) if every proper
holomorphic subbundle $\F$ of ${\mathcal E}$ satisfies $\mu(\F) <
\mu({\mathcal E})$ (respectively $\mu(\F) \leq \mu({\mathcal E})$).
Note that when $n$ and $d$ are coprime, every semistable bundle of
rank $n$ and degree $d$ is stable.

Let $E$ be a fixed complex vector bundle of rank $n$ and degree $d$
over $C$.  The group $Aut_{C}(E)$ of all complex vector bundle
automorphisms of $E$ is called the {\em complexified gauge group} of
$E$ in the Yang-Mills context (see below) and denoted ${\mathscr
G}_{\cplx}$.  Let ${\mathscr C}={\mathscr C}(n,d)$ denote the space
of all holomorphic structures on $E$, and let ${\mathscr C}^{s}$
(respectively ${\mathscr C}^{ss}$) be the subset of ${\mathscr C}$
consisting of stable (respectively semistable) holomorphic
structures on $E$.  Since $C$ is one complex-dimensional, all
almost-complex structures on $E$ are automatically integrable, and
holomorphic structures on $E$ are specified by elements of an
infinite-dimensional complex affine space whose vector space of
translations is isomorphic to the space $\Omega^{0,1}(C,End(E))$.
The complexified gauge group ${\mathscr G}_{\cplx}$ acts on
${\mathscr C}$ by bundle automorphisms, and isomorphism classes of
holomorphic vector bundles on $C$ are in bijection with ${\mathscr
G}_{\cplx}$-orbits in ${\mathscr C}$.  In modern terms, we can
identify the {\em quotient stack} ${\mathscr C}/{\mathscr
G}_{\cplx}$ with the {\em moduli stack} of vector bundles on $C$,
but to obtain a well behaved moduli {\em space} we restrict
ourselves to stable bundles.

Let $\overline{{\mathscr G}_{\cplx}}$ denote the quotient of
${\mathscr G}_{\cplx}$ by its central subgroup, isomorphic to
$\cplx^*$, corresponding to scalar multiples of the identity
automorphism of $E$. This quotient group $\overline{{\mathscr
G}_{\cplx}}$ acts freely on ${\mathscr C}^s$ (which equals
${\mathscr C}^{ss}$ when $n$ and $d$ are coprime) and the quotient
space ${\mathscr C}^s/\overline{{\mathscr G}_{\cplx}}$ can be
naturally identified with $\mnd$, the moduli space of stable
holomorphic bundles on $C$ of rank $n$ and degree $d$.

Atiyah and Bott identify the space ${\mathscr C}$ of holomorphic
structures on $E$ with a space of unitary connections on $E$ in
order to apply Yang-Mills theory to study $\mnd$. A $U(n)$-principal
bundle on $C$ will be called a {\em unitary bundle}; unitary bundles
on $C$ are classified by homotopy classes of maps from $C$ to the
classifying space $BU(n)$.  The set of isomorphism classes of
unitary bundles and the set of isomorphism classes of
$GL_n(\cplx)$-principal bundles are in bijection as the inclusion
$U(n) \hookrightarrow GL_n(\cplx)$ is a homotopy equivalence and
hence induces a homotopy equivalence $BU(n) \longrightarrow BGL_n$.
Thus, a $GL_n(\cplx)$-bundle $E$ admits the structure of a unitary
bundle and all such structures are equivalent up to
$GL_n(\cplx)$-bundle automorphisms.  A Hermitian structure on a
complex vector bundle $E$ is a choice of unitary bundle structure
underlying the given $GL_n(\cplx)$-bundle structure; let us fix such
a Hermitian structure on $E$. The gauge group ${\mathscr G}$ is then
the group of unitary bundle automorphisms of $E$; it is homotopy
equivalent to its complexification which is the group of complex
automorphisms of $E$  we have already denoted by  ${\mathscr
G}_{\cplx}$.  The automorphism group of any unitary bundle contains
a central subgroup isomorphic to $U(1)$.  Let $\overline{{\mathscr
G}}$ denote the quotient of ${\mathscr G}$ by this central subgroup;
then $\overline{{\mathscr G}_{\cplx}}$ is the complexification of
$\overline{{\mathscr G}}$.

There is a canonical affine linear isomorphism of the space
${\mathscr A}$ of unitary connections on the complex vector bundle
$E$ with the space ${\mathscr C}$ of holomorphic structures on $E$
(\cite{AB} p. 570).  The Yang-Mills functional on ${\mathscr A}$
associates with any unitary connection $A$ on $E$ the integral over
$C$ of the norm square of its curvature $F_A$. Atiyah and Bott apply
the ideas of equivariant Morse theory to the Yang-Mills functional
on ${\mathscr A}$ (or equivalently on ${\mathscr C}$).

Recall that the equivariant cohomology $H^*_G(Y)$ of a topological space $Y$ on which a group
$G$ acts can be defined as the ordinary cohomology of the Borel quotient
$EG \times_G Y$ where $EG \to BG$ is the universal $G$-bundle over the classifying
space $BG$ of $G$.
When $G$ acts freely on $Y$ the natural map $EG \times_G Y \to Y/G$
has contractible fibres and induces  an isomorphism $H^*(Y/G) \cong H^*_G(Y)$, while
when $Y$ is contractible the map $EG \times_G Y \to BG$ induces an isomorphism
between $H^*_G(Y)$ and $H^*(BG)$.
Since ${\mathscr C}$ is an affine space it is contractible, and
so
\begin{equation*}
H^*_{{\mathscr G}_{\cplx}}({\mathscr C}) \cong H^*(B{\mathscr G}_{\cplx}).
\end{equation*}
This algebra is easy to describe explicitly: over $\Q$ it is freely
generated as a polynomial algebra tensored with an exterior algebra
by the K\"{u}nneth  components of the equivariant Chern classes of
the universal bundle on ${\mathscr C} \times X$ and has Poincar\'{e}
series
\begin{equation}
\label{eqn:PBG}
\begin{split}
P_t(B{\mathscr G}_{\cplx}) &\stackrel{def}{=} \sum_{i \geq 0} t^i \dim_\Q H^i(B{\mathscr G}_{\cplx},\Q)\\
&= \frac{\prod_{j=1}^n (1 + t^{2j-1})^{2g}}{(1-t^{2n})\prod_{j=1}^{n-1}(1-t^{2j})^2}.
\end{split}
\end{equation}

Since $\overline{{\mathscr G}_{\cplx}}$ acts freely on ${\mathscr
C}^s$, the identification of smooth manifolds $\mnd \cong {\mathscr
C}^{s}/{\mathscr G}_{\cplx} \cong {\mathscr
C}^{s}/\overline{{\mathscr G}_{\cplx}}$ induces the following
isomorphisms on cohomology:
$$
H^*(\mnd) \cong H^*({\mathscr C}^{s}/\overline{{\mathscr
G}_{\cplx}}) \cong H^{*}_{\overline{{\mathscr G}_{\cplx}}}({\mathscr
C}^{s})
$$
and
$$H^{*}_{{{\mathscr G}_{\cplx}}}({\mathscr C}^{s}) \cong H^*(B\cplx^*) \otimes
H^{*}_{\overline{{\mathscr G}_{\cplx}}}({\mathscr C}^{s}).
$$

Atiyah and Bott study the cohomology of the moduli space $\mnd$ when
$n$ and $d$ are coprime (so that ${\mathscr C}^s = {\mathscr
C}^{ss}$) by showing that the restriction map
$H^*_{\overline{{\mathscr G}_{\cplx}}}({\mathscr C}) \longrightarrow
H^{*}_{\overline{{\mathscr G}_{\cplx}}}({\mathscr C}^{ss})$ is
surjective.  In order to do this, they consider the Yang-Mills (or
Atiyah-Bott-Shatz) stratification of ${\mathscr C}$; this is a
stratification of ${\mathscr C}$ by $\overline{{\mathscr
G}_{\cplx}}$-stable submanifolds, which is the \lq Morse
stratification' induced by the Yang-Mills functional but also has a
purely algebraic description as follows.

All holomorphic vector bundles on $C$ are algebraic.  If $C$ is an
algebraic curve over any field $k$, and ${\mathcal E}$ is an
algebraic vector bundle over $C$, then ${\mathcal E}$ has a
canonical {\em Harder-Narasimhan filtration}.  This is an increasing
filtration by algebraic subbundles
\begin{equation}
\label{HNfilt}
0 = F_0({\mathcal E}) \subset F_1({\mathcal E}) \subset \cdots \subset F_r({\mathcal E}) = {\mathcal E}
\end{equation}
uniquely determined by the conditions that $F_j({\mathcal
E})/F_{j-1}({\mathcal E}) = gr_j({\mathcal E})$ is semistable for $1
\leq j \leq r$, of degree $d_j$ and rank $n_j$, say, and that the
sequence of slopes $\mu(gr_j({\mathcal E})) = d_j/n_j$ satisfies
$$
\mu(gr_1({\mathcal E})) > \mu(gr_2({\mathcal E}))  > \cdots > \mu(gr_r({\mathcal E})).
$$
The Harder-Narasimhan type of $E$ is defined to be the decreasing sequence $\underlinemu=
(\mu_1, \ldots ,\mu_n)$ of rational numbers in which $\mu(gr_j({\mathcal E}))
= d_j/n_j$ appears $n_j$ times.
The Yang-Mills stratum
${\mathscr C}_{\mu}$ is then the subset
 of ${\mathscr C}$
consisting of all holomorphic structures of Harder-Narasimhan type $\underlinemu$ on $E$.

Since the Harder-Narasimhan filtration of a holomorphic vector
bundle is canonically defined, the assignments ${\mathcal E} \mapsto
gr_j {\mathcal E}$ induce a map from the space of vector bundles
with a fixed Harder-Narasimhan type $(d_1/n_1, \ldots, d_r/n_r)$ as
above to the product of the spaces of bundles of degree $d_j$ and
rank $n_j$.  If we let ${\mathscr G}_{\cplx}(n_j,d_j)$ denote the
complexified gauge group of the smooth vector bundle underlying
$gr_j {\mathcal E}$, the map just defined (together with the
K\"{u}nneth formula) gives an isomorphism:
\begin{equation}
\label{eqn:strata}
H^*_{{\mathscr G}_{\cplx}}({\mathscr C}_{\underlinemu}) \cong \bigotimes_{j=1}^r H^*_{{\mathscr G}_{\cplx}(n_j,d_j)}({\mathscr C}^{ss}(n_j,d_j)).
\end{equation}
The unique open stratum of the Yang-Mills stratification is ${\mathscr C}^{ss}$, and when $n$ and $d$ are coprime we have
\begin{equation*}
H^*_{{\mathscr G}_{\cplx}}({\mathscr C}^{ss}) \cong H^*(B\cplx^*) \otimes H^*(\mnd).
\end{equation*}

The complexified gauge group ${\mathscr G}_{\cplx}$ acts on
${\mathscr C} = {\mathscr C}(n,d)$ preserving the Yang-Mills
stratification. The Yang-Mills strata may be totally ordered so that
the closure of a stratum ${\mathscr C}_{{\underlinemu}}$ is
contained in the union of lower-dimensional strata ${\mathscr
C}_{{\underlinemu}'}$ with ${\underlinemu}' \geq {\underlinemu}$. We
let $U_{{\underlinemu}}$ denote the open subset of ${\mathscr C}$
obtained by taking the union of all strata ${\mathscr
C}_{{\underlinemu}'}$ for ${\underlinemu}' \leq {\underlinemu}$.  We
can then consider the inclusion of ${\mathscr C}_{\underlinemu}$
into $U_{{\underlinemu}}$ and the associated Thom-Gysin sequence:
$$
\cdots \to H_{{\mathscr G}_{\cplx}}^{j-2c_{{\underlinemu}}}({\mathscr C}_{{\underlinemu}}) \to H^j_{{\mathscr G}_{\cplx}}(U_{{\underlinemu}}) \to
H^j_{{\mathscr G}_{\cplx}}(U_{\underlinemu}-{\mathscr C}_{{\underlinemu}}) \to \cdots
$$
where $c_{\underlinemu}$ is the codimension of the complement of $U_{\underlinemu}$ in ${\mathscr C}_{n,d}$.
The Yang-Mills stratification is equivariantly perfect in the sense that these Thom-Gysin sequences break up into short exact sequences
$$
0 \to H_{{\mathscr G}_{\cplx}}^{j-2c_{\underlinemu}}({\mathscr C}_{{\underlinemu}}) \to H^j_{{\mathscr G}_{\cplx}}(U_{\underlinemu}) \to
H^j_{{\mathscr G}_{\cplx}}(U_{\underlinemu}-{\mathscr C}_{{\underlinemu}}) \to 0.
$$
The integer $c_{\underlinemu}$ can be computed in terms of ranks and
degrees appearing in the Harder-Narasimhan type
${\underlinemu}=(d_1/n_1, \ldots ,d_r/n_r)$:
\begin{equation}
\label{12}
c_{\underlinemu}= \sum_{\ell>j} (n_{\ell}d_{j}-n_{j}d_{\ell}+n_{\ell}n_{j}(g-1)).
\end{equation}
Atiyah and Bott show that the Yang-Mills stratification is
equivariantly perfect by considering the composition of the
Thom-Gysin map $H_{{\mathscr
G}_{\cplx}}^{j-2c_{\underlinemu}}({\mathscr C}_{\underlinemu}) \to
H^j_{{\mathscr G}_{\cplx}}(U_{\underlinemu})$ with restriction to
${\mathscr C}_{\underlinemu}$, which is multiplication by the
equivariant Euler class $e_{\underlinemu}$ of the normal bundle to
${\mathscr C}_{\underlinemu}$ in ${\mathscr C}$. They show that
$e_{\underlinemu}$ is not a zero-divisor in $H_{{\mathscr
G}_{\cplx}}({\mathscr C}_{\underlinemu})$ and deduce that the
Thom-Gysin maps $H_{{\mathscr
G}_{\cplx}}^{j-2c_{\underlinemu}}({\mathscr C}_{\underlinemu}) \to
H^j_{{\mathscr G}_{\cplx}}(U_{\underlinemu})$ are all injective.

Choosing splittings of the short exact sequences above then gives
the following direct sum decomposition (of rational vector spaces):
\begin{equation}
\label{eqn:directsum}
H^j_{{\mathscr G}_{\cplx}}({\mathscr C},\Q) \cong \bigoplus_{i} H_{{\mathscr G}_{\cplx}}^{j-2c_{\underlinemu}}({\mathscr C}_{\underlinemu},\Q).
\end{equation}
One may derive the inductive formulas obtained by Atiyah and Bott in
\cite{AB} for the equivariant Betti numbers of  $\calc^{ss}_{n,d}$
(and hence when $n$ and $d$ are coprime for the Betti numbers of
$\mnd$) by combining equations (\ref{eqn:strata}), (\ref{eqn:PBG})
and (\ref{eqn:directsum}). In terms of the equivariant Poincar\'{e}
series
$$
P_t^{{\mathscr G}_{\cplx}}(X) \stackrel{def}{=} \sum_{i \geq 0} t^i \dim_\Q H^i_{{\mathscr G}_{\cplx}}(X,\Q) $$
they are given by
\begin{equation}
\label{eqn:ptg} P_t^{{\mathscr G}_{\cplx}}({\mathscr C}^{ss}) =
P_t(B{\mathscr G}_{\cplx}) - \sum_{\underlinemu \neq (d/n,\ldots,d/n)} P_t^{{\mathscr G}_{\cplx}}({\mathscr
C}_{\underlinemu})
\end{equation}
where $P_t(B{\mathscr G}_{\cplx})$ is given by (\ref{eqn:PBG}) and
if $\underlinemu = (d_1/n_1, \ldots , d_r/n_r))$ then by
(\ref{eqn:strata}) we have
$$P_t^{{\mathscr G}_{\cplx}}({\mathscr C}_{\underlinemu}) = \prod_{j=1}^r
P_t^{{\mathscr G}_{\cplx}(n_j,d_j)}({\mathscr C}^{ss}(n_j,d_j)).$$

\subsubsection*{The case $n=2$.}
Every line bundle over $C$ is stable, and the moduli spaces ${\mathcal M}(1,d)$
are topologically tori $\cplx^{g}/\Z^{2g}$; we have
$$P_t({\mathcal M}(1,d)) = P_t(B \overline{\mathscr G}_{\cplx}(1,d))=(1+t)^{2g}.$$
When $n=2$ and $d$ is odd, the Harder-Narasimhan filtration of an unstable bundle
${\mathcal E}$ is very simple: it is given by
$$ 0 \subset {\mathcal L} \subset {\mathcal E}$$
where ${\mathcal L}$ is a line subbundle of ${\mathcal E}$ of degree $d_1 > d/2$.
Thus the inductive formula (\ref{eqn:ptg}) takes the form
\begin{equation}
\label{eqn:PBG1}
 \begin{split}
P_t^{{\mathscr G}_{\cplx}(2,d)}({\mathscr C}^{ss}(2,d))   &=
P_t(B{\mathscr G}_{\cplx}(2,d))  \\ & -
\sum_{d_1 > d/2} t^{2(2d_1 - d + g - 1)}
P_t^{{\mathscr G}_{\cplx}(1,d_1)}({\mathscr C}^{ss}(1,d_1))
P_t^{{\mathscr G}_{\cplx}(1,d-d_1)}({\mathscr C}^{ss}(1,d-d_1))
\\
&= \frac{(1 + t)^{2g}(1+t^3)^{2g}}{(1-t^{4})(1-t^{2})^2}
- \sum_{j=0}^{\infty} t^{2(2j+g)} \left( \frac{(1+t)^{2g}}{(1-t^2)} \right)^2.
\end{split}
\end{equation}
It follows that when $d$ is odd
\begin{equation}
\label{eqn:PBG2}
\begin{split}
P_t({\mathcal M}(2,d)) = (1-t^2)
P_t^{{\mathscr G}_{\cplx}(2,d)}({\mathscr C}^{ss}(2,d)) &=
\frac{(1+t)^{2g}}{(1-t^2)} \left( \frac{(1+t^3)^{2g}}{(1-t^4)}
- \frac{t^{2g} (1+t)^{2g}}{(1-t^4)} \right)
 \\
&= \frac{(1+t)^{2g}((1+t^3)^{2g} - t^{2g}(1+t)^{2g})}{(1-t^2)^2(1+t^2)}
\end{split}
\end{equation}
which is a polynomial of degree $6g-6$ in $t$ if $g \geq 2$.

\subsection{Tamagawa numbers and curves over finite fields}
\subsubsection*{The Weil conjectures}
The computation of the Betti numbers of the moduli spaces $\mnd$ of
bundles of coprime rank $n$ and degree $d$ over a nonsingular
complex projective curve $C$ by Atiyah and Bott was preceded by
equivalent inductive formulae presented by Harder, Narasimhan,
Desale and Ramanan \cite{HN,DR}.  These computations utilize the
Weil conjectures, proved by Deligne, and a computation of the
Tamagawa number for $SL_n$, perhaps originally due to Siegel, but
appearing explicitly in work of Weil \cite{Weil}.  Let us briefly
review the setting of the Weil conjectures as it serves to motivate
``motivic" ideas.

Let $\Fq$ be a finite field with $q$ elements, where $q$ is a power
of a prime $p$.  Suppose that $X$ is a finite type scheme over
$\Fq$.  Following Weil, one defines the zeta function of $X$ as the
following formal power series with rational coefficients:
$$
Z_X(t) = \exp(\sum_{r = 1}^{\infty} |X({\mathbb F}_{q^r})|\frac{t^r}{r}).
$$
Thus $Z_X(t)$ is a generating function for the numbers
$|X({\mathbb F}_{q^r})|$
of ${\mathbb F}_{q^r}$-rational points of $X$.

Suppose now that $X$ is a smooth projective variety over $\Fq$ of
dimension $n$.  Weil conjectured that $Z_X(t)$ has the following
three properties, the second and third of which are analogous to
properties of the Riemann zeta function.  First,
$Z_X(t)$ is a rational function of $t$, i.e. it is a quotient of
polynomials with rational coefficients.   Second, if $E$ is the
self-intersection number of the diagonal $\Delta$ of $X \times X$,
then $Z_X(t)$ satisfies a functional equation:
$$
Z_{X}(\frac{1}{q^n t}) = \pm q^{\frac{nE}{2}}t^{E}Z_X(t).
$$
Thirdly, $Z_X(t)$ satisfies an analogue of the Riemann hypothesis.  More precisely, it is possible to write
$$
Z_X(t) = \frac{P_1(t)P_3(t) \cdots P_{2n-1}(t)}{P_0(t)P_2(t) \cdots P_{2n}(t)}
$$
where $P_{0}(t) = 1-t$, $P_{2n}(t) = 1 - q^n t$ and all the other
$P_i(t)$ are polynomials with integer coefficients that can be
written
$$
P_i(t) = \prod (1 - \alpha_{ij}t)
$$
with $\alpha_{ij}$ some collection of algebraic integers of norm
$q^{i/2}$. (The closer analogue of the Riemann zeta function is
$\zeta_X(s)=Z_X(q^{-s})$, whose zeros and poles are on the lines
$\text{Re }s=j/2$ with $j=0,1,\ldots,2n$.) The first two statements
of the Weil conjectures were verified (for arbitrary $n$) first by
Dwork and later by Artin and Grothendieck using Grothendieck's
theory of \'etale cohomology.  The Riemann hypothesis was later
established by Deligne \cite{De1}.

Let $\overline{\Fq}$ denote an algebraic closure of $\Fq$.  Given a
finite type scheme $X$ over $\Fq$, let $\bar{X}$ be the variety over
$\overline{\Fq}$ obtained by base change.  Let $Fr_q$ denote the
geometric Frobenius automorphism of $\bar{X}$ induced by the
automorphism $\alpha \mapsto \alpha^q$ of $\overline{\Fq}$.  Suppose
that $\ell$ is a prime number not equal to $p$, and let
$\Q_{\ell}$ denote the field of $\ell$-adic numbers.  Then Artin and
Grothendieck defined \'etale cohomology groups
$H^i_{\et}(\bar{X},\Q_{\ell})$ satisfying many properties analogous
to usual cohomology groups: these groups are only non-vanishing when
$0 \leq i \leq 2n$, satisfy Poincar\'e duality and have the usual
exact sequences (Mayer-Vietoris, Thom-Gysin, etc.).  Furthermore,
they showed using an $\ell$-adic version of the Lefschetz fixed
point theorem that $P_i(t)$ could be interpreted as the
characteristic polynomial of $Fr_q$ acting on the $\ell$-adic
\'etale cohomology group $H^i_{\et}(\bar{X},\Q_{\ell})$.  Thus
%$\ell$-adic Betti numbers described above, presented as
the degrees of the polynomials $P_i(t)$ are in fact the ranks of these $\ell$-adic
cohomology groups.

Let $L \hookrightarrow \cplx$ be an algebraic number field.  If $X$ is a smooth projective
variety defined over $L$, then one can use the Weil conjectures to
determine the ordinary Betti numbers of $X$ thought of as a complex
variety.  The essential point is that such an $X$ can be treated via base change
both as a variety over a finite field $\Fq$, for appropriate $q$, and
also (given a choice of embedding
$L \hookrightarrow \mathbb{C}$)
as a variety over
$\mathbb{C}$, and furthermore the $l$-adic Betti numbers of the
former agree with the ordinary Betti numbers of the latter. Note
that this is very much an observation about Betti numbers rather than
more refined topological invariants. Indeed, Serre has produced in
\cite{Serre} examples of smooth projective varieties $X$ over $L$
such that simply choosing two different embeddings of $L$ into
$\mathbb{C}$ yields two smooth complex varieties which are
not homeomorphic to each other, although by the Weil conjectures
they must have the same Betti numbers
(more completely, it is known that all such ``conjugate" complex
varieties must have the same ``\'etale homotopy type").

Let us be a little more precise. Suppose that $L$ is an algebraic
number field equipped with a fixed embedding $\phi: L
\hookrightarrow \cplx$ and let ${\mathcal O}_L$ denote the ring of
integers in $L$; we write ${\mathfrak p}$ for a non-zero prime ideal
(necessarily maximal) of $\O_L$ and the quotient $\O_L/{\mathfrak
p}$ is a finite field $\Fq$ (where $q$ is a power of some prime
$p$).  Suppose that $X$ is a smooth projective variety defined over $L$;
we will now consider the $\Fq$-variety that is the ``reduction
modulo ${\mathfrak p}$" of $X$. Choose a projective embedding of
$X$. Upon clearing denominators, which involves only finitely many
prime ideals, we obtain a non-empty open subscheme $B'$ of $\Spec
\O_L$ and an associated morphism of schemes ${\mathfrak X}'
\rightarrow B'$ which is proper (indeed projective) over $B'$ and
restricts to $X$ over the generic point $\Spec L$ of $B'$. Since the
morphism is generically smooth, shrinking $B'$ if necessary we
obtain a non-empty open subscheme $B$ of $\Spec \O_L$ and associated morphism
$\mathfrak{X} \rightarrow B$ which is both smooth and projective and
which restricts to $X$ over the generic point $\Spec L$. Given a
prime ideal $\mathfrak{p}$ whose inclusion in $\O_L$ induces a
morphism $\Spec \O_L/{\mathfrak p} \cong \Spec \Fq \rightarrow B$,
the image is necessarily a closed point of $B$.  The
scheme-theoretic fibre $\mathfrak{X}_{\mathfrak{p}}$ of
$\mathfrak{X} \rightarrow B$ over this closed point is the
``reduction modulo ${\mathfrak p}$" of ${\mathfrak X}$.  Being a
fibre product with $\Spec \Fq$, it is necessarily an $\Fq$-scheme,
and by the construction of $B$ it is both smooth and projective.  Then
the degree of the polynomial $P_i(t)$ defined as above for
${\mathfrak X}_{\mathfrak p}$ is equal to the $i$-th Betti number of
the complex variety $X_{\cplx}$ obtained from $X$ by base change via
the given morphism $\phi^*: \Spec \cplx \longrightarrow \Spec L$; we
will abuse notation and refer to these as Betti numbers of $X$.

Now when $C$ is a smooth projective curve over $\Fq$, geometric
invariant theory (GIT) can be used to define the
moduli space ${\mathcal M}^C(n,d)$ of stable, rank $n$, degree $d$
vector bundles over $C$, at least for large enough $q$
(we will review this construction in section \S 5). 
When $n$ and $d$ are coprime, the resulting moduli space ${\mathcal
M}^C(n,d)$ is again a smooth projective variety over $\Fq$ (or a
suitable finite field extension). 

It is not hard to see that the topology of the moduli space $\mnd$
of stable vector bundles of rank $n$ and degree $d$ over a compact
Riemann surface (or nonsingular complex projective curve) depends
only on the genus of the Riemann surface, not on its complex
structure.  Thus, in order to study the cohomology of $\mnd$, we can
begin with a smooth projective curve $\pi: {\mathfrak C} \longrightarrow B$ over an
open subscheme $B$ of the ring of integers in a number field $L \hookrightarrow \cplx$. Assuming that $\pi$ admits a section, 
Seshadri's extension of GIT to general base schemes (see
\cite{Seshadri}) can be used to construct a moduli space of stable
vector bundles of coprime rank and degree over ${\mathfrak C}$ (see
\cite{Maru} or \cite{Gasb} Th\'eor\`eme 4.3). Given a prime $p$ for
which there exists a morphism $\Spec \Fq \longrightarrow B$ (for $q$
a power of $p$) whose image we denote by ${\mathfrak p}$, base
change produces a scheme isomorphic to the moduli scheme of stable
bundles over a nonsingular projective curve ${\mathfrak
C}_{\mathfrak p}$ over $\Fq$; similarly, using the inclusion of the
generic point $\Spec L \hookrightarrow B$ and the given map $\Spec
\cplx \longrightarrow \Spec L$, base change produces a scheme
isomorphic to the moduli scheme of stable bundles over a nonsingular
projective curve over $\cplx$.  Thus the Weil conjectures allow us
as above to compute the Betti numbers of $ \mnd$ from the numbers of
points defined over ${\mathbb F}_{q^r}$ (for $r=1,2,\ldots$) of
${\mathcal M}^{C}(n,d)$ where $C={\mathfrak C}_{{\mathfrak p}}$ (cf.
\cite{HN} pp. 239-242).

In fact it is technically easier to work with the moduli spaces
${\mathcal M}_{\Lambda}(n,d)$ and ${\mathcal M}^C_{\Lambda}(n,d)$ of
stable vector bundles of rank $n$ and fixed determinant line bundle
$\Lambda$ (of degree $d$); when $n$ and $d$ are coprime, calculating
the Betti numbers of ${\mathcal M}_{\Lambda}(n,d)$ is equivalent to
calculating those of $\mnd$ since the determinant map $\mnd
\rightarrow \mathcal{M}(1,d)$ defined by $E \mapsto \det(E)$ is a
fibration with fibre at $\Lambda \in \mathcal{M}(1,d)$ given by
${\mathcal M}_{\Lambda}(n,d)$, and it induces an isomorphism
$$H^*(\mnd) \cong H^*(M_{\Lambda}(n,d)) \otimes H^*(\mathcal{M}(1,d)).$$
Here ${\mathcal M}(1,d)$ is isomorphic to the Jacobian variety of
$C$, which is an abelian variety of dimension equal to the genus $g$
of $C$.  In this context, counting the number of ${\mathbb
F}_{q^r}$--points of ${\mathcal M}^C_{\Lambda}(n,d)$ amounts to
counting isomorphism classes of stable rank $n$ vector bundles
defined over ${\mathbb F}_{q^r}$ on the curve $C$ with fixed
determinant line bundle $\Lambda$.

\subsubsection*{Tamagawa Numbers}
Let $C$ be a smooth projective curve over a finite field $\Fq$.
Harder, Narasimhan, Desale and Ramanan (see \cite{HN,DR}) showed how
to count the number of ${\mathbb F}_{q^r}$-points of the moduli
space ${\mathcal M}^C_{\Lambda}(n,d)$ of stable rank $n$ degree $d$
vector bundles over $C$ with fixed determinant $\Lambda$ of degree
$d$ on $C$, when $n$ and $d$ are coprime.  In order to do this, they
made use of the Tamagawa number of $SL_n$.  Let us recall the
relationship between the Tamagawa number of $SL_n$ and moduli spaces
of vector bundles.

Suppose $K = \Fq(C)$ is the function field of $C$.  Recall that the
ad\`ele ring ${\mathbb A}_K$ of $K$ is defined as follows.  If $x$
is a closed point of $C$, denote the local ring at $x$ by $\Ox$.
Denote the completion of the local ring at $x$ by $\Ohat$ and the
field of fractions of $\Ohat$ by $\Khat$.  Then $\Ohat$ is a compact
topological ring, and choosing a local parameter determines an
isomorphism to ${\mathbb F}_{q^r}[[t]]$, where $r$ is a strictly
positive integer.  Similarly, $\Khat$ is a locally compact field
isomorphic to ${\mathbb F}_{q^r}((t))$.  For any finite set $S$ of
closed points of $C$, we define ${\mathbb A}_S$ to be the product
$$
{\mathbb A}_S = \prod_{x \in S} \Khat \times \prod_{x \in C - S} \Ohat,
$$
where on the right hand side, the notation $x \in C - S$ should be
read ``$x$ is a closed point of $C - S$."  Observe that each
${\mathbb A}_S$ is a locally compact topological ring.  The sets $S$
are partially ordered by inclusion and we let ${\mathbb A}_K$ be the
locally compact topological ring obtained by taking the colimit of
the rings ${\mathbb A}_S$ as $S$ varies.

Since $SL_n$ is defined over $\Spec \Z$, we can consider the set
$SL_n({\mathbb A}_K)$.  This set can be equipped naturally with the
structure of a locally compact topological group.  Let ${\mathfrak
K}$ be the maximal compact subgroup of $SL_n({\mathbb A}_K)$ which
is the product $\prod_{x} SL_n(\Ohat)$ (again running over closed
points of $C$).  The group $SL_n(K)$ can be viewed as a discrete
subgroup of $SL_n({\mathbb A}_K)$.  Let ${\sf Bun}_{SL_n}(\Fq)$
denote the set of isomorphism classes of vector bundles on $C$,
defined over $\Fq$, with trivial determinant.  The starting point
for the relationship between moduli spaces of bundles over $C$ and
Tamagawa numbers was Weil's construction of a  canonical bijection:
$$
{\mathfrak K} \backslash SL_n({\mathbb A}_K) / SL_n(K) \;\;\;\isomto \;\;\;{\sf Bun}_{SL_n}(\Fq).
$$

Briefly, any vector bundle ${\mathcal E}$ on $C$ is Zariski locally
trivial.  Thus, ${\mathcal E}$ can be trivialized at the generic
point $\eta$ of $C$.  Similarly, ${\mathcal E}$ can be trivialized
over the formal disc $\Spec \Ohat$ for any $x \in C$.  Elements of
$SL_n({\mathbb A}_K)$ can then be identified with collections
$({\mathcal E},\varphi_{\eta}, \setof{\varphi_{x}}_{x \in C})$
consisting of a vector bundle ${\mathcal E}$ over $C$ with trivial
determinant, a trivialization $\varphi_{\eta}$ of ${\mathcal E}$ at
the generic point, and a trivialization $\varphi_{x}$ of ${\mathcal
E}$ over the formal disc at every closed point $x \in C$.  Given an
element $g \in SL_n({\mathbb A}_K)$, the points $x \in C$ at which
$g_x$ does not lie in $SL_n(\Ohat)$ form a closed subscheme $S$ of
$C$.  We can reconstruct the vector bundle ${\mathcal E}$ from $g$
by twisting a trivial bundle on $C - S$ by the $g_x$ for $x \in S$.
For a fixed determinant line bundle $\Lambda \neq {\mathcal O}$ we simply
replace  ${\mathfrak K}$ with a different maximal compact subgroup
of $SL_n({\mathbb A}_K)$ and replace ${\sf Bun}_{SL_n}(\Fq)$ with
the set ${\sf Bun}_{SL_n}^\Lambda(\Fq)$
 of isomorphism classes of vector bundles on $C$, defined over $\Fq$, with fixed determinant
$\Lambda$. One may also consider the double coset space ${\mathfrak
K} \backslash SL_n({\mathbb A}_K) / SL_n(K)$ with $SL_n$ replaced by
other groups.  The corresponding double coset space plays a key
motivating r\^{o}le in the geometric Langlands program (see e.g.
\cite{Gaitsgory,Frenkel}).

As the group $SL_n({\mathbb A}_K)$ is locally compact it possesses a
right invariant Haar measure, which is determined up to scalars. In
fact, given a right invariant, non-zero differential form $\omega$
of top degree on $SL_n$ (thought of as a group scheme over $K$),
there exists a procedure to construct a uniquely determined right
invariant measure $\omega^{\tau}_{{\mathbb A}_K}$ on $SL_n({\mathbb
A}_K)$. This measure induces a measure on the coset space
$SL_n({\mathbb A}_K)/SL_n(K)$ and we set
$$
\tau(SL_n) = \int_{SL_n({\mathbb A_K)/SL_n(K)}} \omega^{\tau}_{{\mathbb A}_K}.
$$
The crucial fact for our discussion is the fact due essentially to
Siegel, but explicitly proved by Weil, that the Tamagawa number
$\tau(SL_n) = 1$ (see \cite{Weil} Theorem 3.3.1).

The differential form $\omega$ induces measures on the compact
groups $SL_n(\Ohat)$ for each $x \in C$, and thus on ${\mathfrak
K}$.   The connection between the Tamagawa number and isomorphism
classes of bundles is provided by Siegel's mass formula:
$$
\tau(SL_n) = vol({\mathfrak K}) \sum_{{\mathcal E} \in {\sf Bun}^\Lambda_{SL_n}(\Fq)} \frac{1}{|Aut({\mathcal E})|}.
$$
If we let $\zeta_C(s) = Z_C(q^{-s})$, the volume of ${\mathfrak K}$
can be computed explicitly in terms of the zeta function of $C$:
$$
vol({\mathfrak K}) = q^{-(n^2-1)(g-1)}\zeta_C(2)^{-1} \cdots \zeta_C(n)^{-1}.
$$
(For generalizations of this result to other groups $G$ --- and to
the case of smooth reductive group schemes over $C$ --- we refer the
reader to the work of Behrend and Dhillon
\cite{BD05b,BD05a}.)

The automorphisms of a stable bundle over $\Fq$ are simply given by
multiplication by nonzero scalars, i.e. elements of the
multiplicative group $\gm(\Fq)$ of $\Fq$.  The set ${\mathfrak K}
\backslash SL_n({\mathbb A}_K)/SL_n(K)\cong {\sf Bun}_{SL_n}^\Lambda
(\Fq) $ can be partitioned by Harder-Narasimhan type into subsets
${\sf Bun}_{SL_n}^{\Lambda,\underlinemu}(\Fq)$.  Via the map
$SL_n({\mathbb A}_K)/SL_n(K) \longrightarrow {\mathfrak K}
\backslash SL_n({\mathbb A}_K)/SL_n(K)$, this induces a
stratification of the coset space $SL_n({\mathbb A}_K)/SL_n(K)$.  It
follows from the Siegel mass formula that when $n$ and $d$ are
coprime the number of isomorphism classes of stable (equivalently
semistable) vector bundles on $C$ with rank $n$ and fixed
determinant $\Lambda$ of degree $d$ is given by
$$
q^{(n^2-1)(g-1)}\zeta_C(2)\cdots\zeta_C(n) - \sum_{\underlinemu \neq
(d/n,\ldots,d/n)} \sum_{{\mathcal E} \in {\sf
Bun}_{SL_n}^{\Lambda,\underlinemu}(\Fq)} \frac{1}{|Aut({\mathcal
E})|},
$$
where as before $\underlinemu$ is the Harder-Narasimhan type of
${\mathcal E}$ determined by the ranks $n_j$ and degrees $d_j$ of
the subquotients of the Harder-Narasimhan filtration of ${\mathcal
E}$.

Summarizing, we get an inductive formula for the sum
$$
\sum_{{\mathcal E} \text{ semistable of rank } n \text{ and degree }
d} \frac{1}{|Aut({\mathcal E})|}$$
where ${\mathcal E}$ runs over
the set of isomorphism classes of semistable vector bundles on $C$,
defined over $\Fq$, of rank $n$ and fixed determinant $\Lambda$ of
degree $d$. When $n$ and $d$ are coprime it enables us to calculate
the Betti numbers of ${\mathcal M}_{\Lambda}(n,d)$ and $\mnd$ via
the Weil conjectures.  This leads to an inductive procedure for
calculating the Betti numbers which is formally the same as that
obtained by Atiyah and Bott via equivariant Morse theory.

\section{Equivariant motivic cohomology}
In the next two sections our aim is to extend the circle of ideas
discussed above to the framework of motives, and motivic homotopy
theory as introduced by Morel and Voevodsky (see \cite{VTriCat,MVW,MV}). In order to do this, we need a version of equivariant
motivic cohomology for a linear algebraic group acting on a
projective variety.  In the subsequent two sections we will apply
this theory to the study of the motivic cohomology of GIT quotients
and of moduli spaces of vector bundles over a nonsingular projective
curve. In this section, we discuss the construction and basic
properties of equivariant motivic cohomology; for more details see
\cite{ADK1} (or its predecessor \cite{EG}
for a definition of equivariant motivic cohomology for smooth varieties
using Bloch's higher Chow groups, and the homological
version for arbitrary varieties).

Throughout this section and the next, $k$ will denote a perfect
field of arbitrary characteristic, and all varieties, groups and
schemes will be assumed to be defined over $k$.  If $G$ is a linear
algebraic group over $k$, a variety (respectively, a projective
variety) $X$ equipped with an algebraic $G$-action will be called
$G$-quasiprojective if it admits an ample $G$-equivariant line
bundle.

\subsection{The Borel construction}
To define an equivariant motivic cohomology theory, we emulate the
Borel construction in topology. (We only discuss the theory for
smooth schemes here; it is possible to develop a good theory for
schemes with mild singularities, e.g., semi-normal schemes.) The
construction of an algebraic model for $BG$ described here is
essentially that of \cite{Totaro,MV}. Let $G$ be a linear algebraic
group over $k$. If $\rho: G \longrightarrow GL(V)$ is any faithful
$k$-rational representation of $V$, we can define a space $EG(\rho)$
as follows. Consider the affine space $V^{\oplus n}$ with its
natural $G$-action.  If $n$ is sufficiently large, then $V^{\oplus
n}$ contains an open subscheme $V_n$ on which $G$ acts freely and
such that the quotient $V_n/G$ exists as a smooth quasiprojective
$k$-variety.  The natural map $V^{\oplus n} \hookrightarrow
V^{\oplus n+1}$ induces a $G$-equivariant map $V_n \hookrightarrow
V_{n+1}$ and hence a morphism $V_n/G \longrightarrow V_{n+1}/G$.  We
then define $BG(\rho)$ to be the ind-scheme $\colim_n V_n/G$ and
$EG(\rho) = \colim_n V_n$.  Similarly, if we let $X$ be a smooth
$G$-quasiprojective variety, then we get an inductive system of
spaces $V_n \times_G X$ and we set $X_G(\rho) = \colim_n V_n
\times_G X$.

Let $\Shv_{Nis}(\Sm/k)$ be the category of Nisnevich sheaves of sets
on $\Sm/k$; for brevity we will denote this category by $\Spc$ and
refer to its objects as spaces.  The assignment $X \mapsto
Hom_{\Sm/k}(\cdot,X)$, sending a scheme to its functor of points,
determines a fully faithful embedding $Sm/k \longrightarrow \Spc$.
This follows from the Yoneda lemma together with the observation
that the Nisnevich topology is sub-canonical (i.e., every
representable presheaf is a sheaf).  Thus, given a smooth scheme
$X$, when we refer to ``the space $X$," we will mean the
corresponding functor.  In a similar way, every ind-scheme can be
viewed as a Nisnevich sheaf of sets and we consider the ind-schemes
$X_G(\rho)$ as Nisnevich sheaves.

The motivic homotopy category $\ho{k}$ (respectively the pointed
motivic homotopy category $\hop{k}$) can be constructed from $\Spc$
(respectively the category $\Spc_{\cdot}$ of pointed spaces) by
localization at an appropriate class of weak equivalences in the
sense of model category theory.  In analogy with classical homotopy
theory, objects of the (pointed) motivic homotopy category are
(pointed) spaces and morphisms are (pointed) ``$\aone$-homotopy
classes of maps," appropriately defined.  We will minimize explicit
definitions involving terms whose definitions are formally analogous
to those from classical homotopy theory.  For example, a space $X$
will be called $\aone$-contractible if it is equivalent to $\Spec k$
in $\ho{k}$.  We remark that a Zariski locally trivial smooth
morphism $f: X \longrightarrow Y$ of smooth schemes with
$\aone$-contractible fibres (e.g., affine space fibres) is an
$\aone$-weak equivalence (see \cite{MV} for a precise definition of
$\aone$-weak equivalence and \cite{ADExotic} for more detailed
discussion of this example).  The space $X_G(\rho)$ gives rise to an
object in $\ho{k}$ or $\hop{k}$.  As one expects, the space
$X_G(\rho)$ viewed as an object of $\hop{k}$ is independent of the
choice of faithful representation $\rho$.

\begin{prop}
\label{prop:borelmodel} For any $\rho$, the space $EG(\rho)$ is
$\aone$-contractible in $\hop{k}$.  For any two faithful
representations $\rho,\rho'$ of $G$, there is a canonical
isomorphism $X_G(\rho) \cong X_G(\rho')$ in $\hop{k}$.
\end{prop}

Henceforth, we write $X_G$ for the object in the motivic homotopy
category defined by $X_G(\rho)$ for {\em any} faithful $k$-rational
representation $\rho$.  The space $X_G$ will be called the {\em
motivic Borel construction} for $G$ acting on $X$.  The proof of the
proposition above is conceptually very simple and involves a ``space
level" version of the Bogomolov double fibration construction (see
\cite{EG} Proposition-Definition 1).  Indeed, for two faithful
representations $\rho,\rho'$, we construct a double inductive system
of spaces $X_{n,n'}$ which map to both $V_n \times_G X$ and $V'_{n'}
\times_G X$, where $V_n \times_G X$ denotes the geometric quotient
variety of $V_n \times X$ by the free action of $G$. Using \cite{MV}
\S4 Prop 2.3, and basic results about commutation of colimits, we
show that as both $n,n' \longrightarrow \infty$, the space
$X_{n,n'}$ becomes weakly equivalent to $X_G(\rho)$ and
$X_G(\rho')$.

Thus, mimicking the definition of ordinary equivariant cohomology,
and given the definition of motivic cohomology, one can make the
following definition of equivariant motivic cohomology.

\begin{defn}
\label{defn:equivmotivic}
The {\em equivariant motivic cohomology} $H^{\bullet,\bullet}_G(X,\Z)$ of a smooth $G$-quasiprojective scheme $X$ is defined by the equality
\begin{equation}
H^{\bullet,\bullet}_G(X,\Z) = H^{\bullet,\bullet}(X_G,\Z)
= H^{\bullet,\bullet}(X_G(\rho),\Z)
\end{equation}
for any faithful representation $\rho$ of $G$.
\end{defn}

Proposition \ref{prop:borelmodel} together with a modification of
\cite{VRed} Prop 6.1 then gives rise to the following result which
asserts that motivic cohomology may be computed as a limit of
motivic cohomologies of approximations to the Borel construction.

\begin{prop}
\label{prop:welldefined}
For any faithful representation $\rho$ of $G$, we have an isomorphism
\begin{equation}
H^{\bullet,\bullet}(X_G(\rho),\Z) \cong lim_n H^{\bullet,\bullet}(V_n \times_G X, \Z).
\end{equation}
\end{prop}

One key ingredient in the proof of this result, which will be
extremely useful when studying moduli spaces of bundles, is the
following standard consequence of the existence of Gysin sequence
regarding excising subvarieties from a smooth variety.  Essentially
by definition, the motivic cohomology groups $H^{p,q}(X,\Z)$ vanish
for $q < 0$.  Basic results about Gysin triangles then give rise to
the following result.

\begin{lem}
\label{lem:excision}
Suppose $X$ is a smooth scheme and $Z$ is a closed subscheme of codimension $d$ in $X$.  Then the restriction map
$$
H^{p,q}(X,\Z) \longrightarrow H^{p,q}(X - Z,\Z)
$$
is an isomorphism whenever $q < d$.  In this situation, we will say
that the restriction map is an isomorphism on motivic cohomology of
weight $q < d$.
\end{lem}

\begin{rem}
\label{rem:spectra} The motivic cohomology of a smooth scheme $X$
can be described as the set of $\aone$-homotopy classes from $X$ to
$K(\Z(q),p)$, where $K(\Z(q),p)$ are motivic Eilenberg-Maclane
spaces.  From this point of view, contravariant functoriality for
equivariant motivic cohomology of arbitrary $G$-equivariant
morphisms of smooth schemes is evident.  Furthermore, for different
``motivic spectra," this approach allows a definition of
``Borel-style" generalized equivariant algebraic cohomology
theories, for example Borel-style equivariant algebraic $K$-theory
or equivariant algebraic cobordism.  In particular, if the motivic
spectrum is a motivic ring spectrum, the associated cohomology
theory has a ring structure; this is known for motivic cohomology.
\end{rem}

\subsection{Basic properties}
It follows immediately from Voevodsky's comparison theorem (see
\cite{VCompare} or \cite{MVW}) and Proposition
\ref{prop:welldefined} that the equivariant motivic cohomology we
have defined is isomorphic to the equivariant higher Chow groups
defined by Edidin-Graham in \cite{EG}.  (Note that
equivariant higher Chow groups constitute a Borel-Moore homology theory,
and hence the isomorphism just mentioned can be viewed as a form of
duality between compactly supported homology and cohomology.)  All
natural properties of motivic cohomology (see \cite{VTriCat} \S 3)
immediately adapt to equivariant motivic cohomology: one has
homotopy invariance for morphisms of schemes which are $\aone$-weak
equivalences, equivariant Mayer-Vietoris sequences, equivariant
Thom-Gysin sequences and projective bundle formulae.  As other
examples of properties of Borel-style equivariant cohomology that
adapt to the motivic setting, we have the following results.

Suppose that $H$ is a closed subgroup of $G$.  If $X$ is an
$H$-quasi-projective scheme, we write $G \times_H X$ for the twisted
product space, i.e., the quotient of $G \times X$ by the action of
$H$ defined by $h \cdot (g,x) = (gh^{-1}, h \cdot x)$. (There seem
to be many different ways to denote this space in the literature. We
follow the standard topological convention and hope the reader does
not confuse this with fiber product constructions.)

\begin{lem}
\label{lem:changeofgroups} Suppose that  $X$ is a smooth
$H$-quasiprojective variety and suppose that  $G$ is a connected
linear algebraic group such that $H \subset G$.  Let $Y$ denote the
twisted product space $G \times_H X$.  Then there is a canonical
isomorphism $X_H \cong Y_G$ in $\hop{k}$.
\end{lem}

\begin{lem}
\label{lem:freeaction} Let $X$ be a smooth $G$-quasiprojective
variety with a scheme-theoretically free $G$-action.  Then the
quotient $X/G$ exists as a smooth scheme and the natural map $X_G
\longrightarrow X/G$ is an isomorphism in $\hop{k}$.
\end{lem}

The first lemma is proved in essentially the same way as Proposition
\ref{prop:borelmodel}.  The second lemma follows from \cite{MV} \S4
Lemma 2.8.

We write $\dmeff{\Q}$ for the derived category of effective rational
motivic complexes over $k$.  There is a covariant functor $\Sm/k
\longrightarrow \dmeff{\Q}$ defined by $X \mapsto {\sf M}(X)$ (see
\cite{MVW} Definition 14.1) where ${\sf M}(X)$ is called the
rational motive of $X$.  The object ${\sf M}(X)$ can be thought of
as analogous to the singular chain complex (with rational
coefficients) of the infinite symmetric product of a topological
space $M$, viewed as an object in the derived category of
$\Q$-vector spaces.  Rational motivic cohomology groups can also be
computed as $Hom_{\dmeff{\Q}}({\sf M}(X),\Q(q)[p])$ where $\Q(q)$ is
a certain complex of sheaves in $\dmeff{\Q}$.  In everything that
follows, we work with motives and motivic cohomology with {\em
rational coefficients}.

In this setting, Lemma \ref{lem:freeaction} can be generalized to varieties with finite quotient singularities.

\begin{prop}
\label{prop:finstabs} Let $X$ be a smooth $G$-quasi-projective
variety such that (i) $G$ acts on $X$ with finite stabilizers and
(ii) a geometric quotient $X/G$ exists as a scheme.  The projection
morphism $X_G \longrightarrow X/G$ then gives an isomorphism in
motivic cohomology
$$
H^{\bullet,\bullet}(X/G,\Q) \isomto H^{\bullet,\bullet}_G(X,\Q).
$$
\end{prop}

We also have a version of the Leray-Hirsch theorem.  Recall
that the topological Leray-Hirsch theorem allows one to compute the
(additive) structure of the cohomology of certain fibre bundles
(where the cohomology of the fibre is a free module over the
cohomology of a point).  Analogously, we have the following result.

\begin{prop}
\label{prop:flagbundle} Let $G$ be a split connected reductive group
over $k$.  Let $P$ be a parabolic subgroup defined over $k$. Suppose
that ${\mathscr P} \longrightarrow X$ is a principal $G$-bundle on a
smooth scheme $X$ which has a $k$-rational point. Then
\begin{itemize}
\item $H^{\bullet,\bullet}(G/P,\Q)$ is free as a module over the motivic cohomology of $\Spec k$, and
\item given elements $c_{\alpha}$ of $H^{\bullet,\bullet}({\mathscr P} \times_G G/P,\Q)$ whose restriction to $G/P$ (the fibre over a fixed $k$-rational point) form a basis for $H^{\bullet,\bullet}(G/P,\Q)$ (such elements always exist), one can construct an isomorphism
of modules over the motivic cohomology of $\Spec k$:
$$
H^{\bullet,\bullet}({\mathscr P} \times_G G/P,\Q) \cong H^{\bullet,\bullet}(X,\Q) \tensor_{H^{\bullet,\bullet}(\Spec k,\Q)} H^{\bullet,\bullet}(G/P,\Q).
$$
\end{itemize}
\end{prop}

Let $T \subset B \subset G$ be a choice of a maximal torus $T$
inside a Borel subgroup $B$ of a split connected reductive group $G$
and suppose that  $X$ is a smooth $G$-quasi-projective scheme.  Fix
a faithful $k$-rational representation $\rho$ of $G$, and consider
the corresponding models $X_G(\rho)$.  Observe that $\rho$ gives
faithful $k$-rational representations of $B$ and $T$ as well (by
restriction) and we obtain maps $X_T(\rho) \longrightarrow X_B(\rho)
\longrightarrow X_G(\rho)$, functorial in $X$.  The morphism
$X_T(\rho) \longrightarrow X_B(\rho)$ can be checked to be Zariski
locally trivial with fibres isomorphic to affine space and is thus
an $\aone$-weak equivalence.  If we assume furthermore that $X$ has
a $k$-rational point, then the fibre of $X_B(\rho) \longrightarrow
X_G(\rho)$ over any such point is isomorphic to $G/B$.  The
identification $X_T \cong (G/T \times X)_G$ induces a natural action
of the Weyl group $W$ on $X_T$.  Tracking the action of the Weyl
group and applying the previous proposition, we obtain the following
result.

\begin{thm}
\label{thm:weylinvariants} Suppose that  $X$ is a smooth
$G$-quasi-projective scheme possessing a $k$-rational point. The
natural map $X_T \longrightarrow X_G$ induces an isomorphism of
rings
\begin{equation}
H^{\bullet,\bullet}_G(X,\Q) \longrightarrow H^{\bullet,\bullet}_T(X,\Q)^W.
\end{equation}
\end{thm}

%bd13Jun: removed because this requires a discussion we have not yet added to ADK, so should be checked first
%\begin{rem}
%It is not required that we take motivic cohomology with rational
%coefficients in the above theorem.  Let $t(G)$ be the torsion index
%of $G$ (see \cite{Totaro}), and suppose that  $R$ is any ring in
%which $t(G)$ is invertible.  The above result can be shown to hold
%for motivic cohomology with coefficients in $R$ (see \cite{ADK1}).
%\end{rem}

\begin{ex}
If $T$ is a split torus, we can fix an isomorphism $T \cong (\gm)^{
n}$ where $n = \operatorname{rk} T$ and $\gm$ is the multiplicative
group of the field $k$.  One can check that $B\gm$ is isomorphic to
an infinite dimensional projective space and that $BT$ is isomorphic
to a product of $\operatorname{rk} T$ copies of an infinite
dimensional projective space.  Thus, using the projective bundle
theorem and the K\"unneth formula, one can show that ${\sf M}(BT)
\cong \bigotimes_{i = 1}^{n} (\oplus_{p \geq 0} \Z(p)[2p])$.
\end{ex}

\section{Perfection and equivariant perfection of stratifications}
\subsection{A schematic review of the cohomology of GIT quotients}
\label{ss:schematicreview} Mumford's geometric invariant theory
\cite{GIT} gives a method to construct and study quotients of
certain reductive group actions on algebraic varieties.  To fix
ideas, consider a complex reductive group $G$, a smooth complex
projective variety $X$ equipped with an algebraic action of $G$ and
a $G$-equivariant very ample line bundle $\L$ on $X$.  In this
situation, Mumford introduced a natural $G$-invariant open subset of
\lq semistable' points $X^{ss} \subset X$ for which a projective
(categorical) quotient variety exists; this quotient will be denoted
$X/\!/G$. In general $X/\!/G$ is not an orbit space for the action
of $G$ on $X^{ss}$; however $X^{ss}$ contains an open subset $X^s$
of \lq stable' points for the linear action such that the image of
the restriction to $X^s$ of the quotient map $X^{ss} \to X/\!/G$ is
an open subset of $X/\!/G$ which can be identified naturally with
$X^s/G$. Mumford showed that stability and semistability of a point
$x \in X$ can be tested via one-parameter subgroups $\lambda: \gm
\longrightarrow G$ where $\gm$ is the multiplicative group of
$\cplx$ (the Hilbert-Mumford criterion); thus there is an
effective way to identify the sets  of stable and  semistable
points.

The complement of $X^{ss}$ is called the set of unstable points and
denoted $X^{us}$. To the choice of linearization, there is a
naturally associated ``instability" stratification of $X^{us}$ (see
\cite{Kempf,Hess1,Kir1}).  By taking $X^{ss}$ to be an open stratum,
this extends to a stratification of $X$.

\begin{rem} Suppose that  $T \subset G$ is a maximal torus. The cocharacter
group $\cochar{T}$ is naturally a $\Z$-module and thus we can form
the tensor product $\cochar{T} \tensor_{\Z} \Q$. The Weyl group $W$
acts on $\cochar{T}$ by conjugation. Technically, to define the
stratification, one specifies a $W$-invariant norm $q$ on
$\cochar{T} \tensor_{\Z} \Q$, but this choice will be unimportant
for our purposes.
\end{rem}

It was shown in \cite{Kir1} that this stratification $\setof{S_\beta
: \beta \in {\mathcal B}}$, where the indexing set ${\mathcal B}$ is
a finite set of co-adjoint orbits of $G$, has the following
properties:

\begin{itemize}
\item[{\bf P})] The stratification is rationally $G$-equivariantly perfect, so that there is a ($\Q$-vector space) isomorphism of equivariant cohomology groups
\begin{equation}
\label{induct}
H^{j}_G(X,\Q) \cong  \bigoplus_{\beta \in \mathcal{B}}H^{j -
2d_\beta}_{G}(S_\beta,\Q)
\end{equation}
for all $j \geq 0$, where $d_\beta$ is the (complex) codimension of $S_\beta$ in $X$.
\item[{\bf S1})] The stratum indexed by $0 \in \mathcal{B}$ coincides with the locus $X^{ss}$ of semistable points of $X$ for the linear $G$-action.
\item[{\bf S2})] If $\beta \in \mathcal{B} \setminus \{ 0 \}$ then there is a nonsingular subvariety $Z_\beta$ of $X$, a reductive subgroup $L_\beta$ of $G$ and a linear action of $L_\beta$ on $Z_\beta$, with corresponding semistable locus denoted $Z_\beta^{ss}$, such that
\begin{equation}
H^{\bullet}_G(S_\beta,\Q) \cong H^{\bullet}_{L_\beta}(Z_\beta^{ss},\Q).
\end{equation}
\end{itemize}

\begin{rem} See Theorem \ref{thm:stratification} for more discussion of the
scheme structure of this stratification.
\end{rem}

This provides us with an inductive procedure for calculating the
$G$-equivariant Betti numbers $\dim H^j_G(X^{ss},\Q)$ of $X^{ss}$.
When $X^{ss} = X^s$ (so that the GIT quotient $X/\!/G$ coincides
with the orbit space $X^{ss}/G$ and $G$ acts with only finite
stabilizers on $X^{ss}$) one observes that
\begin{equation}
H^{\bullet}_G(X^{ss},\Q) \cong  H^{\bullet}(X/\!/G,\Q).
\end{equation}
Moreover the Leray spectral sequence for rational cohomology
associated with the fibration
\begin{equation*}
X \times_G EG \longrightarrow BG
\end{equation*}
degenerates because $X$ is smooth and projective (see \cite{De3}), and one obtains an isomorphism of rational vector spaces
\begin{equation}
H^{\bullet}_G(X,\Q) \cong  H^{\bullet}(X,\Q) \tensor_{\Q} H^\bullet(BG,\Q).
\end{equation}

Thus we obtain a method for calculating the Betti numbers of
$X/\!/G$ in terms of the Betti numbers of X and certain smooth
projective subvarieties which turn up inductively, together with the
classifying spaces of $G$ and certain reductive subgroups of $G$;
indeed using versions of Theorem \ref{thm:weylinvariants} and the
Bialynicki-Birula decomposition (see Theorem \ref{thm:BB}) for
ordinary cohomology, we can reduce this to studying the Weyl group
action on the cohomology of the classifying spaces of a maximal
torus $T$ of $G$ and subtori of $T$, together with components of
their fixed point sets on $X$.

An alternative method for obtaining an equivalent inductive
procedure enabling us to calculate the Betti numbers of $X/\!/G$, at
least when $G$ acts freely on $X^{ss}$, is provided by the Weil
conjectures: the stratification allows us to count the semistable
points of associated varieties defined over finite field as the
total number of points minus the sum over $\beta \in {\mathcal B}$
of the number of points in the stratum labelled by $\beta$ (see
\cite{Kir1}\S 15 for more details).

In the remainder of this section, we will show how to adapt the
stratification just described to study the motivic cohomology of GIT
quotients.

\begin{rem}
These techniques for computing the cohomology of the quotient
variety $X/\!/G$ when $X^{ss}=X^s$ were refined and extended to
cover $X^{ss} \neq X^s$ in a series of papers \cite{Kirdesing,KirIHI,KirIHII,JKloc,JKInt,JKKW}.
\end{rem}

\subsection{The Bialynicki-Birula stratification}
Let $k$ be a perfect field, and let $\gm$ be the multiplicative
group over $k$.  Let $X$ be a smooth $\gm$-projective algebraic
$k$-variety.  The fixed-point locus $X^{\gm}$ of the $\gm$-action of
$X$  is in general disconnected, though smooth, and we denote by
$\{Z_i : i \in I\}$ the set of its connected components.  There is a
stratification of $X$ indexed by $I$ whose properties are summarized
in the following theorem.

\begin{thm}[\cite{BB1},\cite{Hess2}]
\label{thm:BB} Let $\{Z_i : i \in I\}$ be the set of connected
components of $X^{\gm}$.  There is a stratification of $X$ by
$\gm$-stable, smooth, locally closed subvarieties $\{Y_i : i \in I
\}$ together with morphisms $Y_i \longrightarrow Z_i$ for $i \in I$
which are $\gm$-equivariant vector bundles.  The inclusion $Z_i
\hookrightarrow X$ factors through the zero section of the bundle
$Z_i \hookrightarrow Y_i$ and the inclusion $Y_i \hookrightarrow X$
for each $i \in I$.
\end{thm}

As observed by Brosnan (see \cite{Brosnan} and the references
therein), the Bialynicki-Birula decomposition actually gives a
decomposition of the integral motive of $X$.  Thus, one obtains the
following result.

\begin{thm}
Suppose that  $X$ is a smooth $\gm$-projective algebraic variety
over $k$.  With notation as in Theorem \ref{thm:BB}, if $c_i$
denotes the codimension of $Y_i$, then one has a decomposition
$$
H^{\bullet,\bullet}(X,\Z) \cong \bigoplus_{i \in I} H^{\bullet-2c_i,\bullet - c_i}(Z_i,\Z)
$$
of modules over the motivic cohomology of $\Spec k$.
\end{thm}

We can also
use the Bialynicki-Birula decomposition to prove a fixed-point
localization theorem in the setting of $\gm$-equivariant motivic cohomology,
in the same style as \cite{Kir1} (see
\cite{ADK1} for more details).    Furthermore, the decomposition
given above holds in any oriented algebraic cohomology theory (cf.
\cite{NeZa}).

\subsection{The instability stratification}
Let $G$ be a split reductive group over $k$ and let $X$ be a smooth
$G$-projective variety. We fix a very ample $G$-linearized line
bundle $\L$ on $X$ and we write $X^{us} = X^{us}(\L)$ for the
complement of the semistable locus $X^{ss} = X^{ss}(\L)$.  In
addition, we write $X \hookrightarrow {\mathbb P}(V)$ for the
projective embedding determined by $\L$.  Henceforth, we will
suppress $\L$.

Next we need more detail about the
stratification of $X^{us}$ than was discussed in \S
\ref{ss:schematicreview}.

If $T \subset G$ is a split maximal torus of $G$, then we let
$\chara{T}$ denote the character group of $T$ and $\cochar{T}$ the
cocharacter group of $T$.

\begin{rem}
The assumption of splitness here is made for simplicity; one can
prove versions of the theorems below without this assumption in
place.
\end{rem}

\begin{thm}
\label{thm:stratification} There is a natural decomposition of $X$
into $G$-invariant subvarieties $S_{\beta}$ labelled by a finite
subset ${\mathcal B}$ of $\cochar{T} \tensor_{\Z} \Q$ which has the
following properties.  To each $\beta \in {\mathcal B}$ there is a
canonically associated parabolic subgroup $P_\beta \subset G$ with a
$k$-defined Levi subgroup $L_\beta \subset P_\beta$ such that
\begin{itemize}
\item there is a smooth $L_\beta$-stable closed subvariety $Z_\beta \subset X$, which is a component of the fixed point locus of a one-parameter subgroup of $T$ representing $\beta$, and
\item there is a $P_\beta$-stable subvariety $Y_\beta \subset X$ and an $L_\beta$-equivariant surjective morphism $Y_\beta \longrightarrow Z_\beta$, which is Zariski locally trivial with fibres isomorphic to affine spaces.
\end{itemize}
Moreover there is a linearization of the induced $L_\beta$-action on
$Z_\beta$ such that if $Z_\beta^{ss}$ denotes the semistable locus
for this $L_\beta$-action on $Z_\beta$ and $Y_\beta^{ss}$ denotes
the fibre product of $Y_\beta$ and $Z_{\beta}^{ss}$ over $Z_\beta$
then $S_\beta$ is the scheme-theoretic image of $G \times_{P_\beta}
Y_\beta^{ss}$ under the multiplication morphism $m_\beta: G
\times_{P_\beta} Y_\beta^{ss} \longrightarrow X$.

$X$ is the disjoint union of the subvarieties $S_\beta$ and $X^{ss}$
coincides with the subvariety $S_0$ labelled by $0 \in {\mathcal
B}$.  Furthermore, one can totally order ${\mathcal B}$ so that
$\setof{S_\beta:\beta \in {\mathcal B}}$ is a stratification in the
sense that $\bar{S_\beta} \subseteq \cup_{\beta' \geq \beta}
S_\beta$ for each $\beta \in {\mathcal B}$. Finally
\begin{itemize}
\item the morphism $m_\beta: G \times_{P_\beta} Y_\beta^{ss} \longrightarrow S_\beta$ is a finite, $G$-equivariant, birational, surjective morphism, and hence an equivariant resolution of singularities, and
\item if, in addition, the action of $G$ on $X$ has the property that for any $\beta \in {\mathcal B}$ and any point $y \in Z_\beta^{ss}$
$$
Lie(P_\beta) = \setof{\xi \in Lie(G)\; | \; \xi_y  \in T_y Y_\beta},
$$
where $\xi_y \in T_y X$ denotes the tangent vector at $y$ given by
the infinitesimal action of $\xi$ (we will say that the linearized
action is {\em manageable}), then $m_\beta$ is an isomorphism.  This
condition is always satisfied when the characteristic of $k$ is $0$.
\end{itemize}

\end{thm}

Over an arbitrary field of characteristic $0$, the result above was
proved by Hesselink (see \cite{Hess1,Hess2}).  In particular, he
showed that all actions in characteristic zero are manageable.  A
proof of the theorem for varieties over an arbitrary perfect field
may be found in \cite{ADK1}.  Without the manageability hypothesis
counterexamples where $m_\beta$ is not an isomorphism can be given
even  for (Frobenius twisted) $SL_2$-actions on $\pone$ over
algebraically closed fields of characteristic $p > 0$.

\subsection{The motivic cohomology of GIT quotients}
With the notations of the previous section in place, we can state
the main result which is a generalization of the main theorem of
\cite{Kir1}.  Suppose that  $X$ is a smooth $G$-projective variety
over $k$.  Let $\setof{S_\beta:\beta\in {\mathcal B}}$ denote the
instability stratification of $X$.  According to Theorem
\ref{thm:stratification}, the indexing set of the instability
stratification can be totally ordered so  that $\bar{S_\beta}
\subseteq \cup_{\beta' \geq \beta} S_\beta$ for each $\beta \in
{\mathcal B}$; we can choose an indexing function $f: {\mathcal B}
\rightarrow \{ 0,1, \ldots, |{\mathcal B}|-1 \}$ commensurate with
this total order.  For simplicity, let us write $S_i$ in place of
$S_{\beta}$ where $i = f(\beta)$.   Let $U^i$ denote the complement
of $\overline{S_i}$ in $X$ for $0 \leq i \leq |{\mathcal B}|-1$ and
let $U^{|{\mathcal B}|} = X$.  There is thus a diagram
$$
U^i \hookrightarrow U^{i+1} \hookleftarrow S_i
$$
where the left inclusion is a $G$-equivariant open immersion and the
right inclusion is a $G$-equivariant closed immersion and $S_i$ is
the complement of $U^i$ in $U^{i+1}$.  If the linearized action is
manageable, then $S_i$ is in fact smooth and isomorphic to $G
\times_{P_i} Y_i^{ss}$. This means that $(S_i)_G \cong
(Y_i^{ss})_{P_i}$ and hence that
$$H^{\bullet,\bullet}_G(S_i,\Z) \cong H^{\bullet,\bullet}_{P_i}(Y_i^{ss}, \Z)
\cong H^{\bullet,\bullet}_{L_i}(Z_i^{ss},\Z)$$
since $P_i$ is $\aone$-homotopy equivalent to $L_i$ and $Y^{ss}_i$ is $\aone$-homotopy
equivalent to $Z_i^{ss}$.
We can therefore consider the equivariant Thom-Gysin sequence
$$
\cdots \longrightarrow H^{\bullet - 2d_i,\bullet-d_i}_{L_i}(Z_i^{ss},\Q) \longrightarrow H^{\bullet,\bullet}_G(U^{i},\Q) \longrightarrow H^{\bullet,\bullet}_G(U^{i-1},\Q) \longrightarrow \cdots .
$$

 The following result is a statement for motivic cohomology analogous to equivariant perfection.

\begin{thm}
\label{thm:maintheorem} Suppose that $G$ is a split reductive group
over $k$.  Let $X$ be a smooth $G$-projective algebraic $k$-variety
with fixed ample $G$-equivariant line bundle $\L$.  Suppose that the
linearized $G$-action on $X$ is manageable.  Let $\{S_\beta:\beta
\in {\mathcal B}\}$ be the stratification of Theorem
\ref{thm:stratification} with the indexing set ${\mathcal B}$
identified with $\{0,1,\ldots,|{\mathcal B}|-1 \}$ as above.  Then
the Thom-Gysin long exact sequences of the inclusion $S_i
\hookrightarrow U^i$ break up into short exact sequences of the form
$$
0 \longrightarrow H^{\bullet - 2d_i,\bullet-d_i}_{G}(S_i,\Q)  \cong
H^{\bullet - 2d_i,\bullet-d_i}_{L_i}(Z^{ss}_i,\Q) \longrightarrow H^{\bullet,\bullet}_G(U^{i},\Q) \longrightarrow H^{\bullet,\bullet}_G(U^{i-1},\Q) \longrightarrow 0.
$$
\end{thm}

Thus the equivariant cohomology of $X=U^{|{\mathcal B}|}$ can be
``reconstructed" from the equivariant cohomology of $X^{ss}=U^0$ and
the $L_\beta$-equivariant cohomologies of $Z_\beta^{ss}$, which are
inductively of the same form.

The proof of this theorem can be obtained by applying the properties
of motivic cohomology discussed in the previous section to modify
the proof given in \cite{Kir1}: this is due to the power and
``topological" nature of motivic cohomology.  Indeed, the result
follows immediately from a motivic version of the ``Atiyah-Bott
lemma" (see \cite{Kir1} Lemma 2.18 or \cite{AB} Prop 13.4) which
guarantees that a certain equivariant motivic Euler class (see
\cite{VRed} \S4 for the definition) is not a zero divisor in
rational motivic cohomology.  The proof of this lemma involves, as
in the Atiyah-Bott case, a reduction to maximal tori (via Theorem
\ref{thm:weylinvariants}) and an identification of the composite of
the Gysin map with restriction to $Z^{ss}_i$ in the exact sequence
of Theorem \ref{thm:maintheorem} with cupping with the Euler class
of the normal bundle.

Theorem \ref{thm:maintheorem} enables us to compute inductively the
equivariant motivic cohomology of $X^{ss}$, and thus using
Proposition \ref{prop:finstabs} to compute the  motivic cohomology
of the quotient $X/\!/G$ when $X^{ss} = X^s$.

\begin{cor}
\label{cor:cohomquotients} Let a reductive group $G$ act on a smooth
$G$-projective variety $X$ with a fixed $G$-linearized line bundle
$\L$.  Suppose that the conditions of Theorem \ref{thm:maintheorem}
hold and that in addition $X^{ss} = X^s$. Then the inclusion $X^{ss}
\hookrightarrow X$ induces a surjection
\begin{equation} \label{4.4.1}
H^{\bullet,\bullet}_G(X,\Q) \rightarrow H^{\bullet,\bullet}_G(X^{ss},\Q) \cong H^{\bullet,\bullet}(X/\!/G,\Q).
\end{equation}
\end{cor}

\begin{rem}
\label{lowdeg}
When a reductive group $G$ acts linearly on a smooth projective variety $X$ with $X^{ss} \neq X^s$,
 we can still use Theorem \ref{thm:maintheorem} to compute inductively the equivariant
motivic cohomology of $X^{ss}$. Then if $d$ is the codimension of
the complement of $X^s$ in $X^{ss}$ we have
$$H^{i,j}(X^s/G,\Q) \cong
H^{i,j}_G(X^{s},\Q) \cong H^{i,j}_G(X^{ss},\Q)  $$
when $j<d$ by Lemma 3.4.
\end{rem}

%bd13Jun: Removed, because the only current clearly correct argument we have uses Q coefficients
%\begin{rem}
%\label{rem3.2.4} The requirement that we take motivic cohomology
%with rational coefficients in Theorem \ref{thm:maintheorem} and
%Corollary \ref{cor:cohomquotients} can be weakened, though the
%situation is now subtler: it suffices to invert a set of primes
%determined by $G$ and its action on $X$ (see \cite{ADK1} for more
%details).
%\end{rem}

\begin{rem}
One can axiomatize the conditions required to make a version of
Theorem \ref{thm:maintheorem} hold for generalized equivariant
motivic cohomology theories in the sense of Remark
\ref{rem:spectra}.  Closely related results have been obtained by
Chai and Neeman (see \cite{ChaiNeeman}). The conditions are
satisfied by, for example, by
 \'etale cohomology and Betti cohomology (see \cite{ADK1}).
\end{rem}

\begin{rem}
The kernel of the surjection (\ref{4.4.1}) can be studied in
different ways. One way, modelled on the results of \cite{Kir1}, is
to use Theorem \ref{thm:maintheorem}. Another is to relate
intersection theory on $X/\!/G$ to equivariant intersection theory
on $ H^{\bullet,\bullet}_G(X,\Q) $; we shall not pursue this here,
but see, for example, \cite{EG,ESChow,JKInt}.
\end{rem}

\section{The motivic cohomology of moduli spaces of bundles over a curve}
In this section, we show how the discussion of  \S 3 and \S 4
can be brought to bear on the study of motivic cohomology for
the moduli spaces $\mnd$ of stable bundles of coprime rank $n$ and
degree $d$ over a smooth projective curve $C$.  We begin by
recalling the GIT construction of $\mnd$ following Newstead \cite{Newstead}; the
emphasis of our discussion is slightly different from existing
treatments as we aim to keep the analogy between the algebraic and
topological categories at the forefront.

To do this, let $C$ be a smooth projective curve over field $k$
which for simplicity we assume is algebraically closed. (The
condition on $k$ can be weakened, but this complicates discussion of
some of the constructions.) In \S 3, we recalled an algebraic
construction of the classifying space $BGL_n$. This space, uniquely
defined as an object in the $\aone$-homotopy category, was
constructed as a limit of smooth quasi-projective varieties
depending on the choice of a faithful representation $\rho$ of
$GL_n$.  If we take $\rho$ to be the standard $n$-dimensional
representation of $GL_n$, then it is easy to see that we have (as
ind-varieties)
$$
BGL_n(\rho) = \colim_\ell Gr(\ell,n)
$$
where $Gr(\ell,n)$ is the Grassmannian of linear $n$-dimensional
quotients of a fixed $\ell$-dimensional vector space. Henceforth, we
suppress $\rho$ and our main object of study will be the space
$$
Map_d(C,BGL_n) \stackrel{def}{=} \colim_\ell Map_d(C,Gr(\ell,n))
$$
of morphisms of degree $d$ from $C$ to $BGL_n$. Here the spaces
$Map_d(C,Gr(\ell,n))$ are (not necessarily smooth)  varieties
which can be identified with  subschemes of the Hilbert scheme of $C
\times Gr(\ell,n)$.

We can construct the moduli space of semistable
bundles on $C$, at least when $d$ is sufficiently large, as a GIT
quotient of (an open subscheme of) $Map_d(C,Gr(m,n))$ where
$$m=d+n(1-g)$$
with respect to an appropriate
linearization of the induced $GL_m$-action.  The set of semistable points which
arises from this construction will be a smooth quasi-projective
variety and when $n$ and $d$ are coprime, its quotient will be the
smooth projective variety $\mndc$.

The construction of the moduli space $\mndc$ from mapping spaces has
the benefit of being closely related to the original Atiyah-Bott
construction involving the classifying space of the gauge group.
This follows from generalizations of Segal's work in \cite{Segal},
which tell us that the inclusion of $Map_d(C,Gr(\ell,n)$ into the
corresponding space $Map^{sm}_d(C,Gr(\ell,n))$ of smooth maps
 is a
cohomology equivalence up to some degree tending to infinity with
$d$. Very roughly speaking, the algebraic mapping spaces
$Map_d(C,Gr(\ell,n))$, and the limiting space $Map_d(C,BGL_n)$, can
be thought of as ``algebraic approximations to Yang-Mills theory."

In order to study the motivic cohomology of $\mndc$ following the
last two sections, na\"ively, we need to understand three things:
\begin{itemize}
\item the stabilizer groups of the action of $GL_m$ on the open subscheme of semistable points
in $Map_d(C,Gr(m,n))$ where $m=d+n(1-g)$ (in particular, we would like all
semistable points to be stable and hence their stabilizer groups to be
finite),
\item the instability stratification of $Map_d(C,Gr(m,n))$ and the equivariant
motivic cohomologies of the unstable strata, and
\item the $GL_m$-equivariant motivic cohomology of $Map_d(C,Gr(m,n))$.
\end{itemize}
Of course there are problems with each of these points. Firstly, the $GL_m$ action on the subset of semistable
points for the required linearization is not effective.
Regarding
the third point, the space $Map_d(C,Gr(m,n))$ is not a smooth
variety and thus the equivariant motivic cohomology is not defined, at least
given the theory developed in \S 3.  This also leads to
problems with the second point, since while one can define an
instability stratification, it will not necessarily have the
properties advertised in Theorem \ref{thm:stratification} and thus
Theorem \ref{thm:maintheorem} will not apply.

The first difficulty is the easiest to address: the centre $\gm \subset
GL_m$ acts trivially on the space of maps and the quotient
group $PGL_m$ acts freely on the open subscheme of semistable points
with quotient the moduli space $\mnd$.  In particular the stabilizer
group at a semistable point, with respect to $PGL_m$, is trivial.
Therefore, we may apply Lemma \ref{lem:freeaction} and deduce that
the { integral} motivic cohomology of $\mnd$ is isomorphic to the
$PGL_n$-equivariant motivic cohomology of the semistable points in
$Map_d(C,Gr(m,n))$, under the appropriate hypotheses on $m,n$ and
$d$.

%bd13Jun: Added a clarifying sentence about instability strata, in the sense of Section 4, and properness; also a paragraph split
To address the second problem we replace $Map_d(C,Gr(m,n))$ with an
open smooth subscheme $R_{n,d}$ (see \S 5.1) which admits $\mndc$ as
a GIT quotient of its semistable points. (Strictly speaking, we use
the semistable points for a projective closure of $R_{n,d}$.)  Thus
one can define an instability stratification with properties much
the same as those discussed in \S4; however, one must be careful
because $R_{n,d}$ is not proper.

For the third problem, understanding the $GL_m$-equivariant
cohomology of $R_{n,d}$, we recall the construction of an auxiliary
space given by Bifet, Ghione and Letizia
(see \cite{BGL}), which is a scheme-theoretic version of Weil's original 1938
description of bundles on curves.  They introduce ind-schemes of {\em matrix
divisors}, which can be thought of as vector bundles on $C$ equipped
with trivializations at the generic point of $C$. The motivic
cohomology of the relevant ind-scheme of matrix divisors is easy to
compute via the Bialynicki-Birula decomposition (Theorem 4.3 above),
since it is an inductive limit of smooth projective varieties over $k$:
the motivic cohomology of spaces of matrix divisors can be reduced
to studying motivic cohomology of symmetric products of curves (see
\S 5.4 below).

In \S 5.3 we apply the procedure of \S 4 in detail to $R_{n,d}$,
reducing the computation of the equivariant motivic cohomology of
$R_{n,d}^{ss}$ (and equivalently, for $n$ and $d$ coprime, the
motivic cohomology of $\mndc$) to the equivariant motivic cohomology
of $R_{n',d'}$ for various $n'\leq n$ and $d'$.  In \S 5.5 we relate
the equivariant motivic cohomology of $R_{n,d}$ to the motivic
cohomology of a space of matrix divisors, which is itself computed
in \S 5.4. Using this comparison, we describe the motivic cohomology
of ${\mathcal M}^C(n,d)$ in terms of the motivic cohomology of
symmetric powers of the curve $C$.

\subsection{GIT construction of $\mnd$} We will follow Newstead \cite{Newstead}.
Let us fix a smooth algebraic curve $C$ of genus $g \geq 2$ over $k$
an algebraically closed field.  If $\L$ is a fixed degree $\lambda$
line bundle on $C$, then tensoring by $\L$ gives an isomorphism
$$
\varphi_{\L} : \mnd \isomto {\mathcal M}(n,d+n\lambda).
$$
Consequently, in what follows we can assume that $d$ is as large as
we want.

Let ${\mathcal Q}_{m,n}$ denote the universal quotient bundle on
$Gr(m,n)$.  A morphism $f: C \longrightarrow Gr(m,n)$ of degree $d$
determines the rank $n$ and degree $d$ bundle $f^*{\mathcal
Q}_{m,n}$ on $C$.  Furthermore, such a morphism determines a
surjection ${\mathcal O}_C^{\oplus m} \longrightarrow f^*{\mathcal Q}_{m,n}$
via pull-back of the defining surjection ${\mathcal O}_{Gr(m,n)}^{\oplus m}
\longrightarrow {\mathcal Q}_{m,n}$.  Thus, we have defined a map
from the set of degree $d$ morphisms $f: C \longrightarrow Gr(m,n)$
to a set of rank $n$ and degree $d$ bundles over $C$ equipped with a
collection of $m$ global generating sections, such that the
$GL_m$-action on $Map_d(C,Gr(m,n))$ corresponds to changing the
basis of generating sections of $f^*{\mathcal Q}_{m,n}$.

By taking $d$ large enough, we may assume that any semistable bundle
${\mathcal E}$ of degree $d$ and rank $n$ over $C$ has the property
that $H^1(C,{\mathcal E}) = 0$ and ${\mathcal E}$ is generated by
its sections.  Then, by the Riemann-Roch theorem,  $\dim
H^0(C,{\mathcal E}) = d + n(1-g)$, so set $m = d + n(1-g)$ and
define an open subscheme
$$
R_{n,d} \subset Map_d(C,Gr(m,n)),
$$
consisting of maps $f: C \longrightarrow Gr(m,n)$ satisfying the following two conditions:
\begin{itemize}
\item[i)] the natural map $H^0(C,{\mathcal O}_C^{\oplus m}) \longrightarrow H^0(C,f^*{\mathcal Q}_{m,n})$ is an isomorphism;
\item[ii)] $H^1(C,f^*{\mathcal Q}_{m,n}) = 0$.
\end{itemize}
Let $R_{n,d}^{s}$ (respectively $R_{n,d}^{ss}$) denote the subset of
$R_{n,d}$ consisting of maps $f$ such that $f^{*}{\mathcal Q}_{m,n}$
is stable (respectively semistable).

\begin{prop}[\cite{Newstead} \S 5]
For any pair $n,d$ with $d$ sufficiently large, the space $R_{n,d}$ is a smooth,
quasi-projective scheme, on which $GL_m$ acts naturally.  The open
subset $R_{n,d}^s$ (respectively $R_{n,d}^{ss}$) can be realized as
the set of stable (respectively semistable) points for an
appropriate linearization of the induced $PGL_m$-action on a
 projective completion of $R_{n,d}$. The group $PGL_m$
acts freely on $R_{n,d}^s$ and the resulting quotient space
$R_{n,d}^{s}/PGL_m$ is isomorphic to the moduli space ${\mathcal
M}^C(n,d).$
\end{prop}

In establishing this proposition, Newstead shows (see
\cite{Newstead} \S 5 and Remark 6.1) that if $N$ is any sufficiently
large integer then $R_{n,d}$ can be embedded as a nonsingular
quasi-projective subvariety of the product $(Gr(m,n))^N$ via the map
\begin{equation}
\label{emb1}
f \mapsto (f(x_1),...,f(x_N))
\end{equation}
for suitable $x_1,...,x_N \in C$.  Observe that $PGL_m$ acts on
$(Gr(m,n))^{ N}$ diagonally.  If $d$ and $N$ are chosen sufficiently
large, then the locus of stable (respectively semistable) points in
the closure $\overline{R_{n,d}}$ of $R_{n,d}$ in $(Gr(m,n))^{ N}$,
for an appropriate linearization of the $PGL_m$-action,
coincides with the locus in $R_{n,d}$ representing stable
(respectively semistable) bundles.
Thus, the notation $R_{n,d}^s$ and $R_{n,d}^{ss}$ serves a
convenient dual purpose. Moreover, for such $d$ and $N$ we have
$$
{\mathcal M}^C(n,d) \cong {R_{n,d}^{ss}}/\!/PGL_m =  \overline{R_{n,d}}^{ss}/\!/PGL_m.
$$

\begin{rem}
There is a natural evaluation morphism $ev: R_{n,d} \times C
\longrightarrow Gr(m,n)$ which, at the level of points, sends a pair
$(f,x)$ corresponding to a map $f: C \longrightarrow Gr(m,n)$ and a
point $x \in C$ to $f(x)$.  The bundle $ev^*{\mathcal Q}_{m,n}$ has
the universal property that for any map $f \in R_{n,d}$ the
restriction of $ev^*{\mathcal Q}$ to $C$ can be canonically
identified with $f^*{\mathcal Q}$.
\end{rem}

\subsection{Finite-dimensional approximations to Yang-Mills theory}
\label{ss:fdapprox} For this subsection, assume that $k = \cplx$. As
in $\S$2.1 let us fix a complex vector bundle $E$ of rank $n$ and
degree $d$, and denote by ${\mathscr C} = {\mathscr C}(n,d)$
 the space of holomorphic structures on $E$.
 We have already remarked that for each
degree $d$, there is a natural inclusion  $Map_d(C,BGL_n)
\hookrightarrow Map^{sm}_{d}(C,BGL_n)$ of algebraic maps into smooth
maps, where the latter space is homotopy equivalent to the
classifying space of the complexified gauge group $B{\mathscr
G}_{\cplx}$.

There is a natural inclusion of Borel constructions
$(R_{n,d})_{GL_m} = EGL_m \times_{GL_m} R_{n,d} \hookrightarrow
B{\mathscr G}_{\cplx}$ defined as follows. Just as $BGL_m$ is
homotopy equivalent to the infinite Grassmannian $Gr(\infty,m)$, so
we can identify $EGL_m$ with the colimit over $\ell$ of the space of
surjective linear maps from $\cplx^{\ell}$ to ${\cplx}^m$. Given a
morphism $f: C \longrightarrow Gr(m,n)$ representing a point of
$R_{n,d}$, and given a surjective linear map $e:\cplx^\ell \to
\cplx^m$, we can define
$$F(f,e):C \to Gr(\ell,m)$$
to be the morphism taking $x \in C$ to  $f(x)$ thought of as an
$n$-dimensional linear quotient of $\cplx^{\ell}$ rather than of
$\cplx^{m}$ by precomposition with $e$. Letting $\ell \to \infty$ this
construction defines a morphism
\begin{equation} \label{eqn:mapBorel}(R_{n,d})_{GL_m} = EGL_m \times_{GL_m} R_{n,d}
\to Map_d(C,BGL_n) \end{equation} and it is shown in \cite{KirMaps}
(see Corollary 7.4 and Lemma 10.1) that the composition of this morphism with the
inclusion of $Map_d(C,BGL_n)$ in $Map^{sm}_d(C,BGL_n) \simeq B{\mathscr G}_{\cplx}$ induces
isomorphisms in cohomology up to some degree which tends to infinity
as $d$ tends to infinity for fixed $n$.

\begin{rem}
\label{rem:Segal} This result is a limiting case of a generalization of Segal's theorem
\cite{Segal} that the inclusion of the space of holomorphic maps
$Map_d(C,\PP^m)$ of degree $d$ from a compact Riemann surface $C$ to
a projective space $\PP^m$ into the corresponding space
$Map^{sm}_d(C,\PP^m)$ of $C^{\infty}$ maps from $C$ to $\PP^m$
induces isomorphisms in cohomology up to degree $(d-2g)(2m-1)-1$.
Here $Map^{sm}_d(C,\PP^m)$ is an infinite-dimensional space
which is independent of $d$ up to homotopy, whereas $Map_d(C,\PP^m)$
is a finite-dimensional algebraic variety whose dimension tends to
infinity with $d$. Segal's theorem tells us that, from the viewpoint
of cohomology, as $d$ tends to infinity the finite-dimensional varieties $Map_d(C,\PP^m)$ are
giving ever better approximations to the infinite-dimensional
space $Map^{sm}_d(C,\PP^m)$.

Segal's result has been generalized to spaces of maps
from $C$ to, for example, flag manifolds and Grassmannians, and more
recently \cite{BHM01} to maps from $C$ to any compact K\"{a}hler manifold
under a holomorphic action of a connected soluble Lie group $S$ with
an open orbit on which $S$ acts freely.
\end{rem}

The infinite-dimensional affine space ${\mathscr C}$ is
contractible, so the natural map from the Borel construction
${\mathscr C}_{{\mathscr G}_{\cplx}} = E{\mathscr G}_{\cplx}
\times_{{\mathscr G}_{\cplx}} {\mathscr C}$ to $B{\mathscr
G}_{\cplx}$ is a homotopy equivalence. Choosing a section gives us
an inclusion
$$Map_d(C,BGL_n) \hookrightarrow B{\mathscr G}_{\cplx} \hookrightarrow {\mathscr C}_{{\mathscr G}_{\cplx}};$$
we can think of this as given by a $C^{\infty}$ identification of
$f^*{\mathcal Q}_{m,n}$ with our fixed complex bundle $E$ for each
$f \in Map_d(C,BGL_n)$ which provides a
holomorphic structure on $E$. Composing this inclusion with
(\ref{eqn:mapBorel}) gives us an inclusion
 of Borel constructions
\begin{equation}
({R_{n,d}})_{GL_m} \hookrightarrow ({{\mathscr C}})_{{\mathscr G}_{\cplx}}
\end{equation}
which induces isomorphisms on cohomology up to arbitrarily high
degree for $d$ sufficiently large. If $\overline{{\mathscr
G}}_\cplx$ denotes the quotient of ${\mathscr G}_\cplx$ by its
central one-parameter subgroup consisting of automorphisms given by
multiplication by nonzero-scalars, this discussion also shows that
there is an induced map
\begin{equation}
i: ({R_{n,d}})_{PGL_m} \hookrightarrow ({{\mathscr C}})_{\overline{{\mathscr G}}_\cplx},
\end{equation}
again inducing isomorphisms on cohomology up to arbitrarily high
degree for $d$ large enough, which is
 compatible with the formation of quotients in the sense that the diagram
\begin{equation}
\label{eqn:ABB}
\xymatrix{
({R_{n,d}^{ss}})_{PGL_m} \ar[d]\ar[r]^{i} & ({{\mathscr C}^{ss}})_{\overline{{\mathscr G}}_\cplx} \ar[d] \\
R_{n,d}/\!/PGL_m \ar[r]^{\cong} & {\mathscr C}/\!/\overline{{\mathscr G}_{\cplx}}
}
\end{equation}
commutes.  If $n$ and $d$ are coprime, both terms in the bottom row are isomorphic to $\mnd$.

While $R_{n,d}$ is not projective, we can still discuss the
instability stratification associated with the linear action of
$PGL_m$ on $R_{n,d}$. Indeed, one can intersect the strata of the
$PGL_m$-action on the singular space $\overline{R_{n,d}}$ (see \S
5.1) with $R_{n,d}$.  In \S 2 we discussed the Yang-Mills
stratification of ${\mathscr C}$. It is shown in \cite{KirMaps} \S
11 that modulo subsets whose codimension tends to infinity as $d$
tends to infinity, the inclusion of Borel constructions above takes
the instability stratification of $R_{n,d}$ to the Yang-Mills
stratification of ${\mathscr C}_{n,d}$. Moreover, although $R_{n,d}$
is not projective, its instability stratification is equivariantly
perfect at least for cohomology up to some degree which tends to
infinity with $d$, and the equivariant cohomology of its unstable
strata can be described inductively in terms of the equivariant
cohomology of $R^{ss}_{\tilde{n},\tilde{d}}$ for varying $\tilde{n}
< n$ and $\tilde{d}$, again up to some degree which tends to
infinity with $d$. Thus the construction of $\mnd$ as a
finite-dimensional GIT quotient $R_{n,d}/\!/PGL_m$ leads to an
alternative derivation of the inductive formulae for calculating the
Betti numbers of $\mnd$ (for more details see \cite{KirMaps}).

\subsection{Equivariant motivic cohomology of strata}
Assume again that $k$ is an arbitrary algebraically closed field.
Let us now study the action of $PGL_m$ on $R_{n,d}$.  Assume
henceforth that $n$ and $d$ are coprime.  Then $PGL_m$ acts freely
on $R_{n,d}^{ss}$ and by Lemma \ref{lem:freeaction} the projection
map $({R^{ss}_{n,d}})_{PGL_m} \longrightarrow \mndc$ induces
isomorphisms in motivic cohomology:
\begin{equation}
\label{eqn:mot1}
H^{\bullet,\bullet}({\mathcal M}^C(n,d),\Z) \isomto H^{\bullet, \bullet}_{PGL_m}(R_{n,d}^{ss}, \Z).
\end{equation}

For a smooth variety $X$, the motivic cohomology group
$H^{2,1}(X,\Z)$ is canonically isomorphic to the Picard group. Thus,
any line bundle $\L$ over $X$ gives a class $c_{2,1}(\L) \in
H^{2,1}(X,\Z)$ which we refer to as its $(2,1)$-chern class.  In
particular, let $\xi$ denote $c_{2,1}({\mathcal O}_{{\mathbb P}^n}(1)) \in
H^{2,1}({\mathbb P}^n,\Z)$.  One can then compute the motivic
cohomology ring of ${\mathbb P}^n$ to be $H^{\bullet,\bullet}(\Spec
k,\Z)[\xi]/\xi^{n+1}$, and there is an analogous theorem for
projectivized vector bundles. Taking the appropriate limit, it
follows that $H^{\bullet,\bullet}(B\gm,\Z) \cong
H^{\bullet,\bullet}(\Spec k,\Z)[[\xi]]$.

Now, the central one-parameter subgroup $\gm \subset GL_m$ acts trivially on $R_{n,d}^{ss}$ and thus we have an $\aone$-weak equivalence
$$
(R_{n,d}^{ss})_{GL_m} \isomto B\gm \times (R_{n,d}^{ss})_{PGL_ m}.
$$
Applying the projective bundle formula, we obtain an isomorphism of rings:
\begin{equation}
\label{eqn:mot2}
H^{\bullet, \bullet}_{PGL_m}(R_{n,d}^{ss}, \Z)[[\xi]] \isomto H^{\bullet, \bullet}_{GL_m}(R_{n,d}^{ss}, \Z).
\end{equation}

We cannot apply the results of \S 4.4 directly to the instability
stratification of $R_{n,d}$ because $R_{n,d}$ is not projective.
However, motivic cohomology has the property that if $Z$ is a
subvariety of codimension $c$ in a smooth variety $X$ then
$$H^{i,j}(X - Z,\Z) \cong H^{i,j}(X,\Z)$$
for weight $j<c$ (see Lemma \ref{lem:excision}). It therefore
follows from a direct adaptation of the arguments of \cite{KirMaps}
\S\S 8-13 that up to some weight which tends to infinity as $d \to
\infty$ for fixed $n$ the restriction map
\begin{equation}
\label{eqn:mot3}
H^{\bullet, \bullet}_{GL_m}(R_{n,d}, \Z) \longrightarrow H^{\bullet, \bullet}_{GL_m}(R_{n,d}^{ss}, \Z)
\end{equation}
is surjective, and that Theorem \ref{thm:maintheorem} applies to the
instability stratification of $R_{n,d}$. Moreover the
instability stratification $\{ S_\beta : \beta \in {\mathcal B} \}$
of $R_{n,d}$ is determined by Harder-Narasimhan type, modulo
subvarieties whose codimension tends to infinity with $d$, in the
following sense. Given any $M>0$, if $d$ is sufficiently large then
to every Harder-Narasimhan type
$$\mu = (d_1/n_1, \ldots, d_r/n_r)$$
as at (\ref{HNfilt}), such that the codimension $c_\mu$ (given in
Equation \ref{12}) of the corresponding Yang-Mills stratum
${\mathscr C}_\mu$  is at most $M$, we can attach an element
$\beta(\mu)$ of the indexing set ${\mathcal B}$ in such a way that
\begin{itemize}
\item[i)] if $\beta \in {\mathcal B}$ is not of the form $\beta(\mu)$ for some
$\mu$ with $c_\mu \leq M$ then the corresponding stratum
$S_\beta$ of the instability stratification of $R_{n,d}$ has codimension greater than $M$;
\item[ii)] outside a subvariety of codimension at least $M$ in $R_{n,d}$ we have
for any $\mu$ with $c_\mu \leq M$ that $f \in S_{\beta(\mu)}$
if and only if $f^*{\mathcal Q}_{m,n}$ has Harder-Narasimhan type $\mu$;
\item[iii)]$S_{\beta(\mu)}$ has codimension $c_\mu$ in $R_{n,d}$ and
$$S_{\beta(\mu)} \cong GL_m \times_{P_{\beta(\mu)}} Y^{ss}_{\beta(\mu)}$$ where
$P_{\beta(\mu)}$ is a parabolic subgroup of $GL_m$ with Levi subgroup
$$L_{\beta(\mu)} \cong \prod_{j=1}^r GL_{m_j}$$
for $m_j = d_j + n_j(1-g)$ where $\mu = (d_1/n_1, \ldots, d_r/n_r)$, and $Y^{ss}_{\beta(\mu)}$ is smooth and is an
$L_{\beta(\mu)}$-equivariant Zariski locally trivial bundle with fibres isomorphic to affine spaces over
$$Z^{ss}_{\beta(\mu)} \cong \prod_{j=1}^r R^{ss}_{n_j,d_j}$$
modulo subvarieties of codimension
at least $M$.
\end{itemize}

Thus, even though $R_{n,d}$ is not projective, nonetheless its instability stratification
satisfies the analogues for motivic cohomology of the properties P, S1 and S2 described in
\S 4.1.
In particular, taking rational coefficients, the kernel of the surjection (\ref{eqn:mot3}),
for weights at most $M$, is, as a rational vector space, isomorphic to
$$ \prod_{\mu \neq (d/n,\ldots,d/n),c_\mu \leq M} H^{\bullet - 2 c_\umu, \bullet - c_\umu}_{GL_m}(S_{\beta(\umu)},\Q)$$
with
\begin{equation}
\label{eqn:mot4}
H^{\bullet, \bullet}_{GL_m}(S_{\beta(\underlinemu)}, \Q) \cong H^{\bullet, \bullet}_{
\prod_{j=1}^r GL_{d_j + n_j(1-g)}}(\prod_{j=1}^r R^{ss}_{n_j,d_j}, \Q).
\end{equation}

\begin{rem}
Motivic cohomology does not have a K\"unneth formula in the sense of
having a convergent K\"unneth spectral sequence for general smooth
schemes (see \cite{DICell} for more information
about K\"unneth spectral sequences in this context).  Therefore
formula (\ref{eqn:mot4}) is not as simple as its topological
counterpart (2.1.3).
\end{rem}

\begin{rem}
To justify this argument, which involves an application of the
methods of Theorem \ref{thm:maintheorem}, we need to check that the
action of $GL_m$ on $R_{n,d}$ is manageable, in the sense of Theorem
4.4.  For this, let us first describe the tangent space to $R_{n,d}$
at a point.  Let ${\mathcal S}_{m,n}$ denote the universal subbundle
of the trivial rank $m$-bundle  on $Gr(m,n)$; as before ${\mathcal
Q}_{m,n}$ denotes the universal quotient bundle.  Consider a closed
point of $R_{n,d}$ defined by $f:C \to Gr(m,n)$, and observe that
the Zariski tangent space to $R_{n,d}$ at this point is canonically
isomorphic to $H^0(C,f^*({\mathcal S}_{m,n}^\vee \otimes {\mathcal
Q}_{m,n}))$.

Suppose that the Harder-Narasimhan filtration of $f^*{\mathcal Q}_{m,n}$ takes the form:
$$
0={\mathcal E}_0 \subset {\mathcal E}_1 \subset \cdots \subset {\mathcal E}_r = f^*{\mathcal Q}_{m,n},$$
with type specified by
$$\mu = (d_1/n_1, \ldots, d_r/n_r).$$
Then (away from a subvariety of $R_{n,d}$ whose
codimension tends to infinity with $d$) $f$ belongs to $Y^{ss}_{\beta(\mu)}$ if and only
if the linear subspaces
$$
0=H^0(C,{\mathcal E}_0) \subset H^0(C,{\mathcal E}_1) \subset \cdots \subset H^0(C,{\mathcal E}_j) \subset \cdots \subset
H^0(C,{\mathcal E}_r) = k^m
$$
are spanned by the subsets $\{e_1, \ldots, e_{m_j}\}$ of the standard basis $\{e_1, \ldots, e_m\}$
of $k^m$. Moreover $f \in Z^{ss}_{\beta(\mu)}$ if and only if in addition each ${\mathcal E}_j$ is
a direct sum
$${\mathcal E}_j = \F_1 \oplus \cdots \oplus \F_j$$
of semistable bundles $\F_j$ of rank $n_j$ and degree $d_j$ where
$H^0(C,\F_j) \subseteq k^m$ is spanned by the subset $\{ e_{m_1 +
\ldots + m_{j-1} +1}, \ldots, e_{m_1 + \ldots + e_{m_j}}\}$ of the
standard basis of $k^m$. Equivalently $f:C \to Gr(m,n)$ is given by
the composition of a map of the form
$$(f_1,\ldots,f_r):C \to \prod_{j=1}^r Gr(m_j,n_j)$$
(for some $f_j \in R^{ss}_{n_j,d_j}$) with the standard embedding of $\prod_{j=1}^r Gr(m_j,n_j)$ in $Gr(m,n)$, where the restrictions
of the universal subbundle ${\mathcal S}_{m,n}$ and quotient bundle ${\mathcal Q}_{m,n}$
on $Gr(m,n)$ to $\prod_{j=1}^r Gr(m_j,n_j)$ are given by
$$\bigoplus_{j=1}^r {\mathcal S}_{m_j,n_j} \quad \mbox{ and } \quad
\bigoplus_{j=1}^r {\mathcal Q}_{m_j,n_j}.$$
Then $f_j^*({\mathcal Q}_{m_j,n_j}) \cong \F_j$ is semistable of rank $n_j$ and degree $d_j$ and
the tangent space to $R_{n,d}$ at $f$ is
$$H^0(C,f^*({\mathcal S}^\vee_{m,n} \otimes {\mathcal Q}_{m,n})) \cong
\bigoplus_{i,j}H^0(C,f^*_i({\mathcal S}^\vee_{m_i,n_i})\otimes f^*_j({\mathcal Q}_{m_j,n_j})),$$
while the tangent space to $Y^{ss}_{\beta(\mu)}$ at $f$ is
$$\bigoplus_{i \leq j}H^0(C,f^*_i({\mathcal S}^\vee_{m_i,n_i})\otimes f^*_j({\mathcal Q}_{m_j,n_j})).$$
If $\xi \in Lie GL_m = \bigoplus_{i,j}(k^{m_i})^\vee \otimes
k^{m_j}$ has decomposition $\xi = (\xi_{i,j})$ with $\xi_{i,j} \in
(k^{m_i})^\vee \otimes k^{m_j}$ then the infinitesimal action of
$\xi$ at $f$ is given by
$$\xi_f = (\xi_{i,j}^f) \in \bigoplus_{i,j}H^0(C,f^*_i({\mathcal S}^\vee_{m_i,n_i})\otimes f^*_j({\mathcal Q}_{m_j,n_j}))$$
where $\xi_{i,j}^f \in H^0(C,f^*_i({\mathcal S}^\vee_{m_i,n_i})\otimes f^*_j({\mathcal Q}_{m_j,n_j}))$
is the image of $\xi_{i,j}$ under the map
\begin{equation} \label{eqn:dij}
(k^{m_i})^\vee \otimes k^{m_j} \to H^0(C,f^*_i({\mathcal S}^\vee_{m_i,n_i})\otimes f^*_j({\mathcal Q}_{m_j,n_j})) \end{equation}
which comes from the bundle surjection
$$(k^{m_i})^\vee \otimes k^{m_j} \otimes {\mathcal O}_C \to f^*_i({\mathcal S}^\vee_{m_i,n_i})\otimes f^*_j({\mathcal Q}_{m_j,n_j})$$
which factors through $(k^{m_i})^\vee \otimes f^*_j({\mathcal
Q}_{m_j,n_j})$. Since $H^0(C,(k^{m_i})^\vee \otimes f^*_j({\mathcal
Q}_{m_j,n_j}))$ is isomorphic to $(k^{m_i})^\vee \otimes H^0(C,
f^*_j({\mathcal Q}_{m_j,n_j})) \cong (k^{m_i})^\vee \otimes k^{m_j}$
and the bundle $f^*_i({\mathcal Q}^\vee_{m_i,n_i})\otimes
f^*_j({\mathcal Q}_{m_j,n_j}) \cong \F_j^\vee \otimes \F_i$ has no
global sections when $i>j$ (as $\F_i$ and $\F_j$ are both semistable
with $\mu(\F_j) > \mu(\F_i)$), it follows from the long exact
sequence
$$\to H^0(C,f^*_i({\mathcal Q}^\vee_{m_i,n_i})\otimes f^*_j({\mathcal Q}_{m_j,n_j}))
\to H^0(C,(k^{m_i})^\vee\otimes f^*_j({\mathcal Q}_{m_j,n_j})) \to
H^0(C,f^*_i({\mathcal S}^\vee_{m_i,n_i})\otimes f^*_j({\mathcal
Q}_{m_j,n_j})) \to
$$
that the map (\ref{eqn:dij}) is injective when $i>j$. Hence the action of
$GL_m$ on $R_{n,d}$ is manageable, at least away from a subvariety of $R_{n,d}$
whose codimension tends to infinity with $d$.
\end{rem}

\subsection{Matrix divisors and ad\`eles}
In order to complete our inductive procedure for understanding the
motivic cohomology of the moduli spaces ${\mathcal M}^C(n,d)$, what
remains is to compute the equivariant motivic cohomology
$H^{\bullet,\bullet}_{GL_m}(R_{n,d},\Z)$.  Taking $k = \cplx$ and
using ordinary cohomology, a method for doing this was sketched in
\S \ref{ss:fdapprox}, but this argument relied on Atiyah and Bott's
computation of the cohomology of the classifying space of the
complexified gauge group together with a generalization of Segal's
theorem \cite{Segal} and must be replaced by something suitably
algebraic for motivic cohomology.  To achieve this, we will relate
the algebraic approximations to Yang-Mills theory to a dual
description in terms of matrix divisors.

A {\em matrix divisor} of rank $n$ on a smooth projective curve $C$
is a locally free subsheaf of the sheaf $k(C)^{\oplus n}$ (where
$K=k(C)$ is the function field of $C$).  Given a vector bundle $\F$,
observe that specifying an inclusion $\F \hookrightarrow
k(C)^{\oplus n}$ is equivalent to specifying a trivialization of
$\F$ at the generic point.

\begin{rem}
Matrix divisors have the following interpretation.  If $K$ is the
function field $k(C)$ of $C$, we can define ${\mathbb A}_{K}$ in a
manner analogous to that for finite fields (see $\S$2.2 above).  For every $k$-point $x$
of $C$, set $\Ohat$ and $\Khat$ to be the completed local ring and
its field of fractions respectively.  For any finite set $S \subset
C(k)$, define ${\mathbb A}_S = \prod_{x \in S} \Khat \times \prod_{x
\in C(k) - S} \Ohat$.  Then set ${\mathbb A}_K$ to be the colimit
of ${\mathbb A}_S$ as $S$ varies through the partially ordered set
of subsets of $C(k)$.  The set ${\sf Bun}_{GL_n}(k)$ of
isomorphism classes of vector bundles of rank $n$ and degree $d$ on
$C$ can be identified with elements of the double coset space
${\mathfrak K} \backslash GL_n({\mathbb A}_K) / GL_n(K)$, where
${\mathfrak K}$ is a subgroup which when $d=0$ is $ \prod_{x \in
C(k)} \Ohat$.  The elements of $GL_n({\mathbb A}_K)$ can be
identified with the space of collections
$(V,\varphi_\eta,\setof{\varphi_x})$ consisting of a vector bundle
$V$ on $C$, a trivialization at the generic point, and a
trivialization along the formal disc at every point $x \in C(k)$.
Thus matrix divisors can be identified with elements of the coset
space ${\mathfrak K} \backslash GL_n({\mathbb A}_K)$.
\end{rem}

For a fixed effective divisor $D$, we can consider the set of
locally free rank $n$  ${\mathcal O}_C$-submodules $\F$ of $k(C)^{\oplus n}$
which are contained in ${\mathcal O}_C(D)^{\oplus n}$.  The set of such
embeddings in fact forms the set of $k$-points of a disconnected
smooth projective variety whose components are indexed by the degree
$d$ of $\F$. Following \cite{BGL}, we denote these components by
$\divnd(D)$.  If $D \leq D'$, then we have a closed embedding
\begin{equation}
\divnd(D) \hookrightarrow \divnd(D').
\end{equation}
We set $\divnd = \colim_{D} \divnd(D)$; then $\divnd$ is a ind-smooth projective variety.

\subsubsection*{Motivic cohomology of the space of matrix divisors}
The motivic cohomology of the space $\divnd(D)$ can be studied using
the Bialynicki-Birula decomposition associated with a generic
one-parameter subgroup $\gm$ of a maximal torus in $GL_n$ (see
Theorem \ref{thm:BB}).

\begin{lem}
\label{lem:mocosym} Let ${\bf d} = (d_1,\ldots,d_n)$ be a vector of
non-negative integers, let $c_{{\bf d}} = \sum_{1 \leq i \leq n}
(i-1) d_i$, and write $|{\bf d}| = d_1 + \cdots + d_n$.  Then
\begin{equation}
\label{Bifmot}
\begin{split}
H^{\bullet,\bullet}(\divnd(D)) \cong \bigoplus_{|{\bf d}| = {d}} H^{\bullet - 2c_{{\bf d}},\bullet-c_{{\bf d}}}(C^{(d_1)} \times \cdots \times C^{(d_n)})
\end{split}
\end{equation}
where $C^{(j)}$ is the $j$th symmetric power of the curve $C$.
\end{lem}

The proof is essentially that given in \cite{Bifet,BGL}. The
automorphism group of the bundle ${\mathcal O}_C(D)$ is $\gm$.  Thus we
obtain an action of the split maximal torus $T \subset GL_n$ of
diagonal matrices on the sum ${\mathcal O}_C(D)^{\oplus n}$ on $\divnd(D)$,
which can be thought of as a torus in $GL_n(k(C))$.  The components
of the fixed-point set of this torus action correspond to matrix
divisors which are direct sums of line bundles of the form
$$
{\mathcal O}_C(D_1) \oplus \cdots \oplus {\mathcal O}_C(D_n)
$$
and, by taking the cokernels of the inclusions of such bundles into
${\mathcal O}_C(D)^{\oplus n}$, give rise to torsion sheaves on $C$.  Thus we
obtain an identification of the fixed-point loci with products of
Hilbert schemes of points on $C$, or equivalently (since $C$ is
one-dimensional) symmetric powers of $C$.

\begin{rem}
The $\ell$-adic cohomology of $\divnd(D)$ is calculated in
\cite{Bifet} (see also \cite{BGL} Proposition 4.2)  using
 the computation above together with Macdonald's computation of the generating function of the cohomology of symmetric products of a curve.
\end{rem}

\begin{rem}
Note that when $j$ is large enough the symmetric power $C^{(j)}$ of
$C$ is a projective bundle over the abelian variety ${\mathcal
M}^C(1,j)$ which is independent of $j$ up to isomorphism. Thus if
$d_1$ is large enough the product  $C^{(d_1)} \times \cdots \times
C^{(d_n)} $ is a projective bundle over ${\mathcal M}^C(1,0) \times
C^{(d_2)} \times \cdots \times C^{(d_n)}$ and hence its motivic
cohomology is independent of $d_1$ up to some level which tends to
infinity with $d_1$. Note also that $c_{{\bf d}}$ as defined in
Lemma \ref{lem:mocosym} for ${\bf d}=(d_1,\ldots ,d_n)$ is
independent of $d_1$. From this and (\ref{Bifmot}) it follows that
the motivic cohomology of $\divnd(D)$ stabilizes as the degree of
$D$ tends to infinity.
\end{rem}

\subsection{Linking maps to Grassmannians with matrix divisors}
Given any morphism $f: C \longrightarrow Gr(m,n)$ representing a
point of $R_{n,d}$, we obtain a surjective morphism ${\mathcal O}_C^{\oplus m}
\longrightarrow f^*{\mathcal Q}_{m,n}$.  Dualizing this morphism
gives an injective map $f^*{\mathcal Q}_{m,n}^{\vee} \hookrightarrow
{\mathcal O}_C^{\oplus m}$ of ${\mathcal O}_C$-modules.  Choosing a morphism
${\mathcal O}_C^{\oplus m} \longrightarrow {\mathcal O}_C^{\oplus n}$ such that the
composite map $f^*{\mathcal Q}_{m,n}^{\vee} \longrightarrow
{\mathcal O}_C^{\oplus n}$ is injective thus gives rise, by definition, to a
matrix divisor.

Recall that in order to complete our inductive procedure for
studying the motivic cohomology of the moduli spaces $\mndc$ it
remains to compute the equivariant motivic cohomology
$H^{\bullet,\bullet}_{GL_m}(R_{n,d},\Z)$, or equivalently the
motivic cohomology of the Borel construction $(R_{n,d})_{GL_m}$
associated with the action of $GL_m$ on $R_{n,d}$. It follows from
the discussion in \S 5.2 that we can identify $(R_{n,d})_{GL_m}$
with the image of the morphism
$$(R_{n,d})_{GL_m} \to Map_d(C,BGL_n) = \colim_\ell\; Map_d(C,Gr(\ell,n))$$
given at (5.2.1). Thus $(R_{n,d})_{GL_m}$ is the colimit
$$\colim_\ell\; (R_{n,d})^\ell_{GL_m}$$
where $(R_{n,d})^\ell_{GL_m}$ is the open subscheme of $Map_d(C,Gr(\ell,n))$ consisting
of maps $F:C \to Gr(\ell,n)$ satisfying

(i) the map of sections $H^0(\alpha_F):k^\ell = H^0(C,{\mathcal
O}_C^{\oplus \ell}) \to H^0(C,F^*({\mathcal Q}_{\ell,n}))$
associated with the natural map $\alpha_F:{\mathcal O}_C^{\oplus
\ell} \to F^*({\mathcal Q}_{\ell,n})$ is surjective, and

(ii) $H^1(C,F^*({\mathcal Q}_{\ell,n}))=0$.

We will complete our inductive procedure by showing that if $d$ is chosen
suitably and $\ell$ and $D$ are sufficiently large then we have canonical
isomorphisms
$$H^{i,j}((R_{n,d})^\ell_{GL_m},\Z) \cong H^{i,j}(\divt(D),\Z)$$
where $\tilde{d} = n \deg D - d$, for all $j$ up to some weight which can
be taken to be arbitrarily high. We will do this by comparing both
$(R_{n,d})^\ell_{GL_m}$ and $\divt(D)$ with an auxiliary scheme $U^\ell_{n,d}$
which is an open subscheme of the product
$$Map_d(C,Gr(\ell,n)) \times ((k^\ell)^\vee \otimes k^n)$$
where we will interpret elements $\varphi \in k^{\ell n}=
((k^\ell)^\vee \otimes k^n)$ as linear maps from the standard
$\ell$-dimensional vector space $k^\ell$ to the $n$-dimensional
vector space $k^n$.  Given any such linear map $\varphi$  and a
scheme $S$, we let $\varphi_S: k^\ell \tensor {\mathcal O}_S \longrightarrow
k^n \tensor {\mathcal O}_S$ denote the corresponding morphism of trivial
vector bundles over $S$.

\begin{defn}
\label{defn:undl} Let $U^\ell_{n,d}$ be the open subscheme of
$Map_d(C,Gr(\ell,n)) \times ((k^\ell)^\vee \otimes k^n)$ consisting
of pairs $(F,\varphi)$ where $F:C \to Gr(\ell,n)$ and
$\varphi:k^\ell \to k^n$ satisfy (i) and (ii) above and also

(iii) the composition of $\varphi_C:{\mathcal O}_C^{\oplus \ell} \to {\mathcal O}_C^{\oplus n}$ with the
dual of $\alpha_F:{\mathcal O}_C^{\oplus \ell}
\to F^*({\mathcal Q}_{\ell,n})$ is an injective map of ${\mathcal O}_C$-modules
$F^*({\mathcal Q}_{\ell,n}^\vee) \to {\mathcal O}_C^{\oplus n}$.

\end{defn}

Let $\tilde{d} = n \deg D - d$. Then we have morphisms
$$\hat{\theta}: U^\ell_{n,d} \to (R_{n,d})^\ell_{GL_m} \quad \mbox{ and } \quad
\tilde{\theta}: U^\ell_{n,d} \to \divt (D)$$
 such that $\hat{\theta}(F,\varphi)=F$ while
$\tilde{\theta}$ sends $(F,\varphi)$ to the matrix divisor
$$F^*({\mathcal Q}^\vee_{\ell,n}) \otimes {\mathcal O}_C(D) \to {\mathcal O}_C(D)^{\oplus n}$$
obtained by tensoring the injective map of ${\mathcal O}_C$-modules
$F^*({\mathcal Q}_{\ell,n}^\vee) \to {\mathcal O}_C^{\oplus n}$, which is
the composition of $\varphi_C:{\mathcal O}_C^{\oplus \ell} \to {\mathcal O}_C^{\oplus n}$ with the
dual of $\alpha_F:{\mathcal O}_C^{\oplus \ell}
\to F^*({\mathcal Q}_{\ell,n})$,
with the identity on
${\mathcal O}_C(D)$.

The fibre of $\hat{\theta}:U^\ell_{n,d} \to (R_{n,d})^\ell_{GL_m}$
over $F:C \to Gr(\ell,n)$ consists of those $\varphi:k^\ell \to k^n$
such that the composition of $\varphi_C:{\mathcal O}_C^{\oplus \ell} \to
{\mathcal O}_C^{\oplus n}$ with the dual of $\alpha_F:{\mathcal O}_C^{\oplus \ell} \to
F^*({\mathcal Q}_{\ell,n})$ is an injective map of ${\mathcal O}_C$-modules
$F^*({\mathcal Q}_{\ell,n}^\vee) \to {\mathcal O}_C^{\oplus n}$. The
codimension in $(k^\ell)^\vee \otimes k^n$ of the complement of any
such fibre tends to infinity as $\ell \to \infty$
 so we have:

\begin{lem} \label{lem:theta} The morphism
$\hat{\theta}: U^\ell_{n,d} \to (R_{n,d})^\ell_{GL_m}$ induces isomorphisms of motivic
cohomology up to a weight which tends to infinity with $\ell$.
\end{lem}

Suppose now that $\F$ is a locally free rank $n$ ${\mathcal O}_C$-submodule of
${\mathcal O}_C(D)^{\oplus n}$ representing a point of $\divt(D)$ and
satisfying the two conditions:

\begin{itemize}
\item[i)] $E^{\vee} \tensor {\mathcal O}_C(D)$ is generated by sections, and
\item[ii)] $H^1(C,E^{\vee} \tensor {\mathcal O}_C(D))$ vanishes.
\end{itemize}

Then the fibre of $\tilde{\theta}:U^\ell_{n,d} \to \divt(D)$ at this point consists of
all pairs $(\varphi, \psi)$ where $\varphi^\vee:k^{\oplus n} \to k^{\oplus \ell}$ is
injective and $\psi:k^\ell \to H^0(C,\F^\vee(D))$ is surjective and
$\psi \circ \varphi:k^{\oplus n} = H^0(C,{\mathcal O}_C^{\oplus n})$ is the map on sections induced by
the dual of $E(-D) \hookrightarrow {\mathcal O}_C^{\oplus n}$. This fibre is therefore
an affine bundle over the space of injective linear maps from $k^n$ to $k^\ell$, which
is an open subscheme of $(k^n)^\vee \otimes k^\ell$ whose complement has codimension tending to
infinity with $\ell$.

Note that the complement in $\divt(D)$ of the open subscheme where these conditions
(i) and (ii) above are satisfied has codimension tending to infinity with $\deg(D)$.
Putting these together yields

\begin{lem}
\label{lem:tildetheta} If $\ell$ and $\deg(D)$ are sufficiently large then
$\tilde{\theta}: U^\ell_{n,d} \to \divt (D)$ induces isomorphisms of motivic
cohomology up to arbitrarily high weights.
\end{lem}

Combining Lemmas \ref{lem:theta} and \ref{lem:tildetheta} we obtain:

\begin{cor}
Given any $M>0$, if $d$ is sufficiently large and if $\ell$ and $\deg(D)$ are large
enough (depending on $d$), we have canonical isomorphisms
 in motivic cohomology
\begin{equation}
\label{eqn:mot5}
H^{i,j}(\divt(D);\Z) \cong H^{i,j}_{GL_m}({R}_{n,d};\Z)
\end{equation}
for $j<M$.
\end{cor}

\begin{rem}
We can think of $(R_{n,d})^\ell_{GL_m}$ mapping to the moduli stack
${\mathscr M}^C(n,d)$ of bundles over $C$ of rank $n$ and degree
$d$, by associating with $F:C \to Gr(\ell,n)$ the bundle
$F^*({\mathcal Q}_{\ell,n})$, with \lq fibre' $EGL_m/\gm = {\mathbb
P}^{\infty}$. Similarly we can think of $\divt(D)$ mapping to the
moduli stack ${\mathscr M}^C(n,\tilde{d})$ where $\tilde{d} = n
\deg(D) - d$ by associating with the matrix divisor $\F
\hookrightarrow {\mathcal O}_C(D)^{\oplus n}$ the bundle $E^\vee
\otimes {\mathcal O}(D)$. The compositions of these with
$\hat{\theta}:U_{n,d}^\ell \to (R_{n,d})^\ell_{GL_m}$ and
$\tilde{\theta}:U_{n,d}^\ell \to \divt(D)$ respectively agree modulo
the isomorphism from ${\mathscr M}^C(n,d)$ to ${\mathscr
M}^C(n,n\deg D - d)$ given by $\F \to \F^\vee \otimes {\mathcal
O}(D)$.
\end{rem}

\begin{rem}
Associating with a matrix divisor $\F \hookrightarrow k(C)^{\oplus
n}$ the bundle $\F$ and forgetting the embedding of $\F$ in
$k(C)^{\oplus n}$ determines a map from the space of matrix divisors
on $C$ to the moduli stack of rank $n$ bundles on $C$; we refer to
this map as the Abel-Jacobi map $\theta$.  When $n$ and $d$ are
coprime, it is shown in \cite{BGL} that if we denote by
$(\divnd)^{ss}$ the space consisting of matrix divisors whose
underlying locally free sheaves are semistable, the Abel-Jacobi map
restricts to a morphism $\theta:  (\divnd)^{ss} \longrightarrow
\mndc$. The isomorphism (\ref{eqn:mot5}) is essentially equivalent
to Theorem 4.5 of \cite{Dh06}, which says that the Abel-Jacobi map
from the ind-variety $\divnd$ to the moduli stack of rank $n$ and
degree $d$ vector bundles on $C$ is a quasi-isomorphism (i.e.
induces isomorphisms on cohomology groups).

The space $(\divnd)^{ss}$ can be identified as an open subset of a
certain projective bundle over $\mndc$.  More precisely,
$(\divnd)^{ss}$ is constructed as an ind-scheme by defining an
inductive system of schemes $\divnd(D)^{ss}$ indexed by effective
divisors $D$ on $C$.  The fibre of $\divnd(D)^{ss}$ at a semistable
vector bundle ${\mathcal E}$ can be identified with the subset of
${\mathbb P}(H^0(C,{\mathcal E}^{\vee} \tensor {\mathcal O}_C(D)^{\oplus n}))$
corresponding to injective ${\mathcal O}_C$-module maps ${\mathcal E}
\hookrightarrow {\mathcal O}_C(D)^{\oplus n}$.  By \cite{BGL} Lemma 8.2, the
codimension of the complement of the fibre of $\divnd(D)^{ss}$ over
${\mathcal E}$ in ${\mathbb P}(H^0(C,{\mathcal E}^{\vee} \tensor
{\mathcal O}_C(D)^{\oplus n}))$ tends to infinity as $deg(D) \longrightarrow
\infty$.  It then follows from the projective bundle theorem for
motivic cohomology that the motivic cohomology of $(\divnd)^{ss}$ is
given by
$$H^{\bullet,\bullet}((\divnd)^{ss}) \cong H^{\bullet,\bullet}({\mathcal M}^C(n,d))[[\xi]]$$
where $\xi$ corresponds to the (2,1) Chern class of the projective
bundle described above over $\divnd(D)^{ss}$ for suitably large $D$.

We can stratify $\divnd(D)$ according to the Harder-Narasimhan type
(\ref{HNfilt}) of the underlying bundle $\F$; this is how Bifet,
Ghione and Letizia obtained their inductive formulae for calculating
the Betti numbers of $\mnd$. Equivalently the ind-variety $\divnd$
can be stratified by Harder-Narasimhan type; the resulting strata
are called Shatz strata in \cite{BGL}.  The Shatz stratification is
perfect and the cohomology of the strata can be described
inductively in a manner identical to our discussion for $R_{n,d}$
but in terms of products of spaces of the form
$Div_{C/k}^{n_j,d_j}(D)^{ss}$ with $0< n_j < n$. If $D \leq D'$, the
closed immersion $\divnd(D) \hookrightarrow \divnd(D')$ is a
stratified morphism, and the limit of these finite dimensional
approximating strata as $\deg(D) \to \infty$ agrees with the Shatz
stratification.

Building on the work of Bifet, Ghione and Letizia, del~Ba{\~n}o \cite{dBmotive} showed, using
the same Shatz stratification of $\divnd$ but now working over a
characteristic 0 not necessarily algebraically closed field, that
one can in the same manner compute the rational Chow motive of the
space of matrix divisors, and thence the virtual Chow motive of the moduli
space of stable bundles on a curve $C$.  To do this he needed to
understand the Chow motive of symmetric products of $C$; he
therefore established a motivic version of MacDonald's formulae for
symmetric products of $C$ \cite[Section 3.3]{dBmotive}. He also needed to show that the Bialynicki-Birula
stratification of a smooth projective variety induces a direct sum
decomposition of the corresponding rational Chow motive; this result
was fundamental to Brosnan's more general theorems on
Bialynicki-Birula decompositions \cite{Brosnan}. He thereby showed that the Chow motive of $\mnd$ lies in the category generated by the motive of the curve and provided a formula for the ``virtual motive" of $\mndc$, with  a closed form expression for the motivic Poincar\'{e} polynomial.
\end{rem}

When $n$ and $d$ are coprime,
 combining (\ref{eqn:mot1}), (\ref{eqn:mot2}), (\ref{eqn:mot3}), (\ref{eqn:mot4}), (\ref{eqn:mot5}) and (\ref{Bifmot})
with the application of Theorem \ref{thm:maintheorem} to the instability stratification
of $R_{n,d}$ as in \S 5.3 leads
to an inductive method for calculating the motivic cohomology groups
$H^{\bullet,\bullet}(\mnd,\Q)$ in terms of the motivic cohomology of products of Jacobians of $C$ and the closely related Hilbert schemes of points on $C$.

%\begin{rem} Some additional care shows that in fact we do not need rational coefficients here (cf. Remark \ref{rem3.2.4} above and \cite{AB}).
%\end{rem}

\begin{rem}
\label{parabolic}
Similar arguments apply to moduli spaces of parabolic bundles over $C$ (cf. \cite{BR,Hollapoincarepoly,HLWeightedframes,MS,Ni,Ni2,NiQPE}.
\end{rem}

\begin{rem}
When $n$ and $d$ have a common factor then (\ref{eqn:mot1}) is no longer valid, but
(\ref{eqn:mot2}), (\ref{eqn:mot3}), (\ref{eqn:mot4}), (\ref{eqn:mot5}) and (\ref{Bifmot}) still provide
an inductive procedure for calculating the equivariant motivic cohomology groups
$H^{\bullet,\bullet}_{GL_m}(R^{ss},\Q)$ of $R^{ss}_{n,d}$ or equivalently the groups $H^{\bullet,\bullet}_{PGL_m}(R^{ss},\Q)$. The codimension of the complement of $R^s$ in $R^{ss}$ is at least
$(g-1)(n-1)$ (cf. e.g. [Kir86b] \S 3 or \cite{Dh06} Cor 5.5) and hence if ${\mathcal M}^{s}(n,d)$ denotes the moduli space of stable bundles we have
$$
H^{i,j}_{PGL_m}(R^{ss},\Q) \cong H^{i,j}_{PGL_m}(R^{s},\Q) \cong H^{i,j}(\calm^s(n,d),\Q)
$$
for $j<(g-1)(n-1)$ by Remark 4.7, giving us the motivic cohomology of
the moduli space $\calm^s(n,d)$ in low degrees. We can also use the
method described in $\S$4.4 above to study the motivic cohomology of a
partial desingularization $\tilde{{\mathcal M}}(n,d)$ of $\mnd$ when $k$ has characteristic zero
(cf. \cite{Kir5}).
\end{rem}

\section{Conclusion}

We have described three equivalent inductive procedures leading to
calculation of the Betti numbers of the moduli spaces $\mnd$ when
$n$ and $d$ are coprime.  The three procedures all rely on
stratifications linked to the Harder-Narasimhan type of a
holomorphic vector bundle, but they differ technically.

(i) The approach of Harder and Narasimhan \cite{HN} and Desale and Ramanan \cite{DR}
via the Weil conjectures uses Tamagawa measures and reduces
to the fact that the Tamagawa number $\tau_{SL_n}$ of $SL_n$ is $1$.

(ii) The approach of Atiyah and Bott \cite{AB} via Yang-Mills theory uses
equivariant Morse theory and reduces to a simple description of the
cohomology of the classifying space of the gauge group.

(iii) The moduli space $\mnd$ can be expressed as a GIT quotient
$R_{n,d}/\!/PGL_m$ and both versions (using equivariant Morse theory and
counting points of an associated variety over finite fields) of the
procedure described in $\S$4.1 for calculating Betti
numbers of GIT quotients $X/\!/G$ can be applied (though extra care is
needed since $R_{n,d}$ is not projective); the $GL_m$-equivariant cohomology
of $R_{n,d}$ can be expressed in terms of the cohomology of symmetric products of the curve
$C$, by relating the Borel construction $(R_{n,d})_{GL_m}$ to a space of matrix divisors.

We have seen in $\S$5 that the third of these approaches can be made to work
for motivic cohomology; it also provides a link between the other
two approaches.  To link the third approach with the approach of
Atiyah and Bott via Yang-Mills theory, we consider the inclusion
$$EGL_m  \times_{GL_m} R_{n,d} \longrightarrow Map_d(C, BSL_n)
\longrightarrow Map^{sm}_d(C, BSL_n) = B{\mathcal G}_{\cplx}$$
defined as at (5.2.1). By a generalization of Segal's theorem on the
topology of spaces of rational maps \cite{Segal} (see Remark
\ref{rem:Segal}) this induces isomorphisms
$$H^j(B {\mathcal{G}}_{\cplx}) \cong H^j_{GL_m}(R_{n,d})$$ for $d$
sufficiently large, which gives us a direct link between approaches (ii) and (iii) (see
\cite{KirMaps} for more details).

To link the third approach with the first using Tamagawa measures,
we consider the variety $U^\ell_{n,d}$ as in Definition
\ref{defn:undl}. The argument given in \S 5.5, using the maps
$\theta: U^\ell_{n,d} \to (R_{n,d})^\ell_{GL_m}$ and
$\tilde{\theta}: U^\ell_{n,d} \to \divt(D)$ to the space of matrix
divisors $\divt (D)$ to show that the $GL_m$-equivariant motivic
cohomology of $R_{n,d}$ is isomorphic to the motivic cohomology of
$\divt(D)$ (up to some weight which can be made arbitrarily large),
also allows us to relate the number of points in corresponding
varieties defined over finite fields. Interpreting matrix divisors
in terms of ad\`eles as in Remark 5.6 then gives us a direct link
between the \lq counting points' version of approach (iii) (see \S
4.1) and the Tamagawa measures used in approach (i).  It also
provides an alternative method for proving that the Tamagawa number
of $SL_n$ is $1$ by using the Bialynicki-Birula stratification of
$\divt(D)$ and counting points on symmetric products of associated
curves over finite fields
 (cf. \cite{BD05a}).

In conclusion, both the arithmetic approach (i) and the Yang-Mills approach (ii)
using equivariant cohomology describe
the moduli space $\mnd$ in terms of an infinite-dimensional quotient
construction; these two infinite-dimensional constructions, though very different,
can be approximated by finite-dimensional quotient constructions which are very closely
related to each other, and in this finite-dimensional setting the Weil conjectures
provide the required link between the arithmetic and the equivariant cohomological
points of view.

%\bibliographystyle{alpha}
%\bibliography{TamagawaFinal}

\begin{thebibliography}{JKKW06}

\bibitem[AB83]{AB}
M.~F. Atiyah and R.~Bott.
\newblock The {Y}ang-{M}ills equations over {R}iemann surfaces.
\newblock {\em Philos. Trans. Roy. Soc. London Ser. A}, 308(1505):523--615,
  1983.

\bibitem[AD07]{ADExotic}
A.~Asok and B.~Doran.
\newblock On unipotent quotients and some {{${\mathbb A}^1$}}-contractible
  smooth schemes.
\newblock {\em Int. Math. Res. Pap.}, 5:1--51, 2007.

\bibitem[ADK]{ADK1}
A.~Asok, B.~Doran, and F.C. Kirwan.
\newblock Equivariant motivic cohomology and quotients.
\newblock {\em In preparation.}

\bibitem[AMM98]{AMM98}
A.~Alekseev, A.~Malkin, and E.~Meinrenken.
\newblock Lie group valued moment maps.
\newblock {\em J. Differential Geom.}, 48(3):445--495, 1998.

\bibitem[AMW00]{AMW00}
A.~Alekseev, E.~Meinrenken, and C.~Woodward.
\newblock Group-valued equivariant localization.
\newblock {\em Invent. Math.}, 140(2):327--350, 2000.

\bibitem[AMW01]{AMW01}
A.~Alekseev, E.~Meinrenken, and C.~Woodward.
\newblock The {V}erlinde formulas as fixed point formulas.
\newblock {\em J. Symplectic Geom.}, 1(1):1--46, 2001.

\bibitem[AMW02]{AMW02}
A.~Alekseev, E.~Meinrenken, and C.~Woodward.
\newblock Duistermaat-{H}eckman measures and moduli spaces of flat bundles over
  surfaces.
\newblock {\em Geom. Funct. Anal.}, 12(1):1--31, 2002.

\bibitem[Ati78]{mfa}
M.~F. Atiyah.
\newblock The unity of mathematics.
\newblock {\em Bull. London Math. Soc.}, 10(1):69--76, 1978.

\bibitem[BB73]{BB1}
A.~Bia{\l}ynicki-Birula.
\newblock Some theorems on actions of algebraic groups.
\newblock {\em Ann. of Math. (2)}, 98:480--497, 1973.

\bibitem[BDa]{BD05b}
K.~Behrend and A.~Dhillon.
\newblock {Connected components of moduli stacks of torsors via Tamagawa
  numbers}.

\bibitem[BDb]{BD05a}
K.~Behrend and A.~Dhillon.
\newblock {On the Motive of the Stack of Bundles}.
\newblock arXiv:math.AG/0512640.

\bibitem[BGL94]{BGL}
E.~Bifet, F.~Ghione, and M.~Letizia.
\newblock On the {A}bel-{J}acobi map for divisors of higher rank on a curve.
\newblock {\em Math. Ann.}, 299(4):641--672, 1994.

\bibitem[BHM01]{BHM01}
C.~P. Boyer, J.~C. Hurtubise, and R.~J. Milgram.
\newblock Stability theorems for spaces of rational curves.
\newblock {\em Internat. J. Math.}, 12(2):223--262, 2001.

\bibitem[Bif89]{Bifet}
E.~Bifet.
\newblock Sur les points fixes du sch\'ema {${\mathrm{ Quot}}_{{\mathscr
  O}^{\oplus r}_X/X/k}$} sous l'action du tore {${\mathbb G}^r_{m,k}$}.
\newblock {\em C. R. Acad. Sci. Paris S\'er. I Math.}, 309(9):609--612, 1989.

\bibitem[Blo86]{Bloch}
S.~Bloch.
\newblock Algebraic cycles and higher {$K$}-theory.
\newblock {\em Adv. in Math.}, 61(3):267--304, 1986.

\bibitem[BR96]{BR}
I.~Biswas and N.~Raghavendra.
\newblock Canonical generators of the cohomology of moduli of parabolic bundles
  on curves.
\newblock {\em Math. Ann.}, 306(1):1--14, 1996.

\bibitem[Bro05]{Brosnan}
P.~Brosnan.
\newblock On motivic decompositions arising from the method of {B}ia\l
  ynicki-{B}irula.
\newblock {\em Invent. Math.}, 161(1):91--111, 2005.
\newblock arXiv:math.AG/0407305.

\bibitem[CN98]{ChaiNeeman}
C-L. Chai and A.~Neeman.
\newblock The naturality in {K}irwan's decomposition.
\newblock {\em Manuscripta Math.}, 97(4):429--434, 1998.

\bibitem[dB01]{dBmotive}
S.~del Ba{\~n}o.
\newblock On the {C}how motive of some moduli spaces.
\newblock {\em J. Reine Angew. Math.}, 532:105--132, 2001.

\bibitem[dB02]{dBrk}
S.~del Ba{\~n}o.
\newblock On the motive of moduli spaces of rank two vector bundles over a
  curve.
\newblock {\em Compositio Math.}, 131(1):1--30, 2002.

\bibitem[Del74]{De3}
P.~Deligne.
\newblock Th\'eorie de {H}odge. {III}.
\newblock {\em Inst. Hautes \'Etudes Sci. Publ. Math.}, (44):5--77, 1974.

\bibitem[Del80]{De1}
P.~Deligne.
\newblock La conjecture de {W}eil. {II}.
\newblock {\em Inst. Hautes \'Etudes Sci. Publ. Math.}, (52):137--252, 1980.

\bibitem[Dhi06]{Dh06}
A.~Dhillon.
\newblock On the cohomology of moduli of vector bundles and the {T}amagawa
  number of {${\rm SL}\sb n$}.
\newblock {\em Canad. J. Math.}, 58(5):1000--1025, 2006.

\bibitem[DI05]{DICell}
D.~Dugger and D.~C. Isaksen.
\newblock Motivic cell structures.
\newblock {\em Algebr. Geom. Topol.}, 5:615--652 (electronic), 2005.

\bibitem[DR75]{DR}
U.~V. Desale and S.~Ramanan.
\newblock Poincar\'e polynomials of the variety of stable bundles.
\newblock {\em Math. Ann.}, 216(3):233--244, 1975.

\bibitem[EG98]{EG}
D.~Edidin and W.~Graham.
\newblock Equivariant intersection theory.
\newblock {\em Invent. Math.}, 131(3):595--634, 1998.

\bibitem[ES89]{ESChow}
G.~Ellingsrud and S.~A. Str{\o}mme.
\newblock On the {C}how ring of a geometric quotient.
\newblock {\em Ann. of Math. (2)}, 130(1):159--187, 1989.

\bibitem[Fre]{Frenkel}
E.~Frenkel.
\newblock {Lectures on the Langlands Program and Conformal Field Theory}.
\newblock arXiv:hep-th/0512172.

\bibitem[Gai03]{Gaitsgory}
D.~Gaitsgory.
\newblock Informal introduction to geometric {L}anglands.
\newblock In {\em An introduction to the Langlands program (Jerusalem, 2001)},
  pages 269--281. Birkh\"auser Boston, Boston, MA, 2003.

\bibitem[Gas97]{Gasb}
C.~Gasbarri.
\newblock Hauteurs canoniques sur l'espace de modules des fibrés stables sur
  une courbe algébrique.
\newblock {\em Bull. Soc. Math. France}, 125(4):457--491, 1997.

\bibitem[Hes78]{Hess1}
W.~H. Hesselink.
\newblock Uniform instability in reductive groups.
\newblock {\em J. Reine Angew. Math.}, 303/304:74--96, 1978.

\bibitem[Hes81]{Hess2}
W.~H. Hesselink.
\newblock Concentration under actions of algebraic groups.
\newblock In {\em Paul Dubreil and Marie-Paule Malliavin Algebra Seminar, 33rd
  Year (Paris, 1980)}, volume 867 of {\em Lecture Notes in Math.}, pages
  55--89. Springer, Berlin, 1981.

\bibitem[HJ94]{HJ94}
J.~Huebschmann and L.~C. Jeffrey.
\newblock Group cohomology construction of symplectic forms on certain moduli
  spaces.
\newblock {\em Internat. Math. Res. Notices}, (6):245 ff., approx.\ 5 pp.\
  (electronic), 1994.

\bibitem[HJ00]{HLWeightedframes}
J.~C. Hurtubise and L.~C. Jeffrey.
\newblock Representations with weighted frames and framed parabolic bundles.
\newblock {\em Canad. J. Math.}, 52(6):1235--1268, 2000.

\bibitem[HN75]{HN}
G.~Harder and M.~S. Narasimhan.
\newblock On the cohomology groups of moduli spaces of vector bundles on
  curves.
\newblock {\em Math. Ann.}, 212:215--248, 1974/75.

\bibitem[Hol00]{Hollapoincarepoly}
Y.~I. Holla.
\newblock Poincar\'e polynomial of the moduli spaces of parabolic bundles.
\newblock {\em Proc. Indian Acad. Sci. Math. Sci.}, 110(3):233--261, 2000.

\bibitem[Jef94]{J}
L.~C. Jeffrey.
\newblock Extended moduli spaces of flat connections on {R}iemann surfaces.
\newblock {\em Math. Ann.}, 298(4):667--692, 1994.

\bibitem[JK95]{JKloc}
L.~C. Jeffrey and F.~C. Kirwan.
\newblock Localization for nonabelian group actions.
\newblock {\em Topology}, 34(2):291--327, 1995.

\bibitem[JK98]{JKInt}
L.~C. Jeffrey and F.~C. Kirwan.
\newblock Intersection theory on moduli spaces of holomorphic bundles of
  arbitrary rank on a {R}iemann surface.
\newblock {\em Ann. of Math. (2)}, 148(1):109--196, 1998.

\bibitem[JKKW03]{JKKW}
L.~C. Jeffrey, Y-H. Kiem, F.~C. Kirwan, and J.~Woolf.
\newblock Cohomology pairings on singular quotients in geometric invariant
  theory.
\newblock {\em Transform. Groups}, 8(3):217--259, 2003.

\bibitem[JKKW06]{JKKW2}
L.~C. Jeffrey, Y-H. Kiem, F.~C. Kirwan, and J.~Woolf.
\newblock Intersection pairings on singular moduli spaces of bundles over a
  {R}iemann surface and their partial desingularisations.
\newblock {\em Transform. Groups}, 11(3):439--494, 2006.

\bibitem[JW94]{JW94}
L.~C. Jeffrey and J.~Weitsman.
\newblock Torus actions, moment maps, and the symplectic geometry of the moduli
  space of flat connections on a two-manifold.
\newblock In {\em Mathematical aspects of conformal and topological field
  theories and quantum groups (South Hadley, MA, 1992)}, volume 175 of {\em
  Contemp. Math.}, pages 149--159. Amer. Math. Soc., Providence, RI, 1994.

\bibitem[Kem78]{Kempf}
G.~R. Kempf.
\newblock Instability in invariant theory.
\newblock {\em Ann. of Math. (2)}, 108(2):299--316, 1978.

\bibitem[Kir84]{Kir1}
F.~C. Kirwan.
\newblock {\em Cohomology of quotients in symplectic and algebraic geometry},
  volume~31 of {\em Mathematical Notes}.
\newblock Princeton University Press, Princeton, NJ, 1984.

\bibitem[Kir85]{Kirdesing}
F.~C. Kirwan.
\newblock Partial desingularisations of quotients of nonsingular varieties and
  their {B}etti numbers.
\newblock {\em Ann. of Math. (2)}, 122(1):41--85, 1985.

\bibitem[Kir86a]{KirMaps}
F.~Kirwan.
\newblock On spaces of maps from {R}iemann surfaces to {G}rassmannians and
  applications to the cohomology of moduli of vector bundles.
\newblock {\em Ark. Mat.}, 24(2):221--275, 1986.

\bibitem[Kir86b]{Kir5}
F.~Kirwan.
\newblock On the homology of compactifications of moduli spaces of vector
  bundles over a {R}iemann surface.
\newblock {\em Proc. London Math. Soc. (3)}, 53(2):237--266, 1986.

\bibitem[Kir86c]{KirIHI}
F.~Kirwan.
\newblock Rational intersection cohomology of quotient varieties.
\newblock {\em Invent. Math.}, 86(3):471--505, 1986.

\bibitem[Kir87]{KirIHII}
F.~Kirwan.
\newblock Rational intersection cohomology of quotient varieties. {II}.
\newblock {\em Invent. Math.}, 90(1):153--167, 1987.

\bibitem[Mac62]{MacDonald}
I.~G. Macdonald.
\newblock The {P}oincar\'e polynomial of a symmetric product.
\newblock {\em Proc. Cambridge Philos. Soc.}, 58:563--568, 1962.

\bibitem[Mar78]{Maru}
M.~Maruyama.
\newblock Moduli of stable sheaves. ii.
\newblock {\em J. Math. Kyoto Univ.}, 18(3):557--614, 1978.

\bibitem[Mei05]{Mei05}
E.~Meinrenken.
\newblock Witten's formulas for intersection pairings on moduli spaces of flat
  {$G$}-bundles.
\newblock {\em Adv. Math.}, 197(1):140--197, 2005.

\bibitem[MFK94]{GIT}
D.~Mumford, J.~Fogarty, and F.~Kirwan.
\newblock {\em Geometric invariant theory}, volume~34 of {\em Ergebnisse der
  Mathematik und ihrer Grenzgebiete (2) [Results in Mathematics and Related
  Areas (2)]}.
\newblock Springer-Verlag, Berlin, third edition, 1994.

\bibitem[MS80]{MS}
V.~B. Mehta and C.~S. Seshadri.
\newblock Moduli of vector bundles on curves with parabolic structures.
\newblock {\em Math. Ann.}, 248(3):205--239, 1980.

\bibitem[MV99]{MV}
F.~Morel and V.~Voevodsky.
\newblock {${\mathbb A}^1$}-homotopy theory of schemes.
\newblock {\em Inst. Hautes \'Etudes Sci. Publ. Math.}, 90:45--143 (2001),
  1999.

\bibitem[MVW06]{MVW}
C.~Mazza, V.~Voevodsky, and C.~Weibel.
\newblock {\em Lecture Notes on Motivic Cohomology}, volume~2 of {\em Clay
  Mathematics Monographs}.
\newblock American Mathematical Society, Providence, RI, Clay Mathematics
  Institute, Cambridge, MA, 2006.

\bibitem[MW99]{MW99}
E.~Meinrenken and C.~Woodward.
\newblock Moduli spaces of flat connections on {$2$}-manifolds, cobordism, and
  {W}itten's volume formulas.
\newblock In {\em Advances in geometry}, volume 172 of {\em Progr. Math.},
  pages 271--295. Birkh\"auser Boston, Boston, MA, 1999.

\bibitem[New78]{Newstead}
P.~E. Newstead.
\newblock {\em Introduction to moduli problems and orbit spaces}, volume~51 of
  {\em Tata Institute of Fundamental Research Lectures on Mathematics and
  Physics}.
\newblock Tata Institute of Fundamental Research, Bombay, 1978.

\bibitem[Nit86]{Ni}
N.~Nitsure.
\newblock Cohomology of the moduli of parabolic vector bundles.
\newblock {\em Proc. Indian Acad. Sci. Math. Sci.}, 95(1):61--77, 1986.

\bibitem[Nit96]{Ni2}
N.~Nitsure.
\newblock Quasi-parabolic {S}iegel formula.
\newblock {\em Proc. Indian Acad. Sci. Math. Sci.}, 106(2):133--137, 1996.

\bibitem[Nit97]{NiQPE}
N.~Nitsure.
\newblock Erratum: ``{Q}uasi-parabolic {S}iegel formula''.
\newblock {\em Proc. Indian Acad. Sci. Math. Sci.}, 107(2):221--222, 1997.

\bibitem[NZ06]{NeZa}
A.~Nenashev and K.~Zainoulline.
\newblock Oriented cohomology and motivic decompositions of relative cellular
  spaces.
\newblock {\em J. Pure Appl. Algebra}, 205(2):323--340, 2006.

\bibitem[Seg79]{Segal}
G.~Segal.
\newblock The topology of spaces of rational functions.
\newblock {\em Acta Math.}, 143(1-2):39--72, 1979.

\bibitem[Ser64]{Serre}
J.-P. Serre.
\newblock Exemples de vari{\'e}t{\'e}s projectives conjugu{\'e}es non
  hom{\'e}omorphes.
\newblock {\em C. R. Acad. Sci. Paris}, 258:4194--4196, 1964.

\bibitem[Ses77]{Seshadri}
C.~S. Seshadri.
\newblock Geometric reductivity over arbitrary base.
\newblock {\em Advances in Math.}, 26(3):225--274, 1977.

\bibitem[Tel98]{Tel98}
C.~Teleman.
\newblock Borel-{W}eil-{B}ott theory on the moduli stack of {$G$}-bundles over
  a curve.
\newblock {\em Invent. Math.}, 134(1):1--57, 1998.

\bibitem[Tot99]{Totaro}
B.~Totaro.
\newblock The {C}how ring of a classifying space.
\newblock In {\em Algebraic $K$-theory (Seattle, WA, 1997)}, volume~67 of {\em
  Proc. Sympos. Pure Math.}, pages 249--281. Amer. Math. Soc., Providence, RI,
  1999.

\bibitem[Voe00]{VTriCat}
V.~Voevodsky.
\newblock Triangulated categories of motives over a field.
\newblock In {\em Cycles, transfers, and motivic homology theories}, volume 143
  of {\em Ann. of Math. Stud.}, pages 188--238. Princeton Univ. Press,
  Princeton, NJ, 2000.

\bibitem[Voe02]{VCompare}
V.~Voevodsky.
\newblock Motivic cohomology groups are isomorphic to higher {C}how groups in
  any characteristic.
\newblock {\em Int. Math. Res. Not.}, (7):351--355, 2002.

\bibitem[Voe03]{VRed}
V.~Voevodsky.
\newblock Reduced power operations in motivic cohomology.
\newblock {\em Publ. Math. Inst. Hautes \'Etudes Sci.}, 98:1--57, 2003.

\bibitem[Wei82]{Weil}
A.~Weil.
\newblock {\em Adeles and algebraic groups}, volume~23 of {\em Progress in
  Mathematics}.
\newblock Birkh\"auser Boston, Mass., 1982.
\newblock With appendices by M. Demazure and Takashi Ono.

\bibitem[Wit92]{Witten}
E.~Witten.
\newblock Two-dimensional gauge theories revisited.
\newblock {\em J. Geom. Phys.}, 9(4):303--368, 1992.

\end{thebibliography}

{{ A. Asok},
University of Washington, Seattle
USA, {asok@math.washington.edu}}

{ { B. Doran},
Institute for Advanced Study, Princeton 
USA, {doranb@math.ias.edu}}

{{ F. Kirwan},
The Mathematical Institute, Oxford
UK, {kirwan@maths.ox.ac.uk}}

\end{document}